\documentclass[a4paper,12pt]{amsart}

\usepackage{amsthm}
\usepackage{amsmath}
\usepackage{amsfonts}
\usepackage{amssymb}
\usepackage{times}
\usepackage{graphicx}

\theoremstyle{definition}
\newtheorem{defn}{\indent\bf Definition}
\newtheorem{rem}[defn]{\indent\bf Remark}
\newtheorem{ex}[defn]{\indent\bf Example}
\theoremstyle{plain}
\newtheorem{lemma}[defn]{\indent\bf Lemma}
\newtheorem{prop}[defn]{\indent\bf Proposition}
\newtheorem{thm}[defn]{\indent\bf Theorem}
\newtheorem{cor}[defn]{\indent\bf Corollary}
\newtheorem*{cor*}{\indent\bf Corollary}

\def\tilde{\widetilde}
\def\a{\alpha}
\def\T{\theta}
\def\PSL{\mathop{\rm PSL}\nolimits}
\def\SL{\mathop{\rm SL}\nolimits}
\def\PGL{\mathop{\rm PGL}}

\def\card{\mathop{\rm card}}
\def\d{{\rm d}}
\def\Z{\mathbb Z}
\def\R{\mathbb R}

\def\cZ{{\mathcal Z}}

\def\Q{\mathbb Q}
\def\N{\mathbb N}

\def\L{{\mathcal L}}

\def\H{{\mathcal H}}

\def\S{{\mathcal S}}
\def\Y{{\mathcal Y}}
\def\O{{\mathcal O}}

\def\B{{\mathcal B}}

\def\cC{{\mathcal C}}

\def\cN{{\mathcal N}}

\def\cX{{\mathcal X}}
\def\cQ{{\mathcal Q}}

\begin{document}

\title[Representations of the Roe algebra of a discrete group and symmetries]
{Unitary representations of the Roe algebra of a discrete group and symmetries}
\author[Florin R\u adulescu]{Florin R\u adulescu${}^*$
\\ \\ 
Dipartimento di Matematica\\ Universit\` a degli Studi di Roma ``Tor Vergata''}

\maketitle 

\thispagestyle{empty}

\renewcommand{\thefootnote}{}
\footnotetext{${}^*$ Member of the Institute of  Mathematics ``S. Stoilow" of the Romanian Academy}
\footnotetext{${}^*$
Supported in part by PRIN-MIUR, and by a grant of the Romanian National Authority for Scientific Research, project number PN-II-ID-PCE-2012-4-0201
}
\footnotetext{${}^*$ Visiting Professor, Department of Mathematics, University of Copenhagen}\begin{abstract}

Let $\Gamma$ be a discrete countable group. Consider the crossed product C$^\ast$-algebra $\mathfrak{R}(\Gamma) = C^{\ast}(\Gamma \rtimes l^{\infty}(\Gamma))$.
Let $G$ be a larger  discrete group, containing $\Gamma$ as an almost normal subgroup.
Consequently $G$ acts by partial isomorphisms on $G$ and hence on $\mathfrak {R}(\Gamma)$.
Let $\mathfrak{R}_G(\Gamma)$ be the crossed product $C^{\ast}$ - algebra $C^{\ast}(G \times (\mathfrak{R}(\Gamma))$. The C$^\ast$-algebra  $\mathfrak{R}_G(\Gamma)$ has a  natural representation into $\B(\ell ^2(\Gamma))$ and hence also admits a representation $\Pi_{\mathcal{Q}}$ into the Calkin algebra $\mathcal{Q}(\ell ^2(\Gamma))$.

Let $G\rtimes \Gamma=\Gamma\times \Gamma^{\rm op} $ and  assume that $\Gamma$ is exact. Assume that the non-trivial conjugation orbits under the action of $\Gamma$, having non amenable stabilizers, are separated, in a suitable chosen  profinite topology, from the identity element in $\Gamma$.  We also assume natural amenability conditions on the dynamics of the action of $\Gamma\times \Gamma^{\rm op}$ on   cosets  of amenable subgroups. Then $\Pi_{\mathcal Q}$ factorises to a representation of $C^{\ast}_{\rm red}(G \rtimes \mathfrak{R}(\Gamma))$.   In particular the groups  $\mathop{SL}_3(\Z)$, $\mathop{\rm PGL}_2(\mathbb Z[\frac{1}{p}])$ have the Akemann-Ostrand property.

This implies, using the solidity property of Ozawa ([Oz]), that  the  group von Neumann algebras,
$\mathcal L(\mathop{SL}_3(\Z))$ and $\mathcal L(\mathop{SL}_n(\Z))$, $n\geq 4$, are non-isomorphic.



\end{abstract}

\vskip0.5cm

\section{Introduction}\label{intro}

Let $G$ be a discrete, countable group.
Assume that $\Gamma \subseteq G$ is almost normal subgroup. Let $G$ act by partial isomorphisms on $\Gamma$. If $\Gamma = G$, then this action is simply the inner action of $\Gamma$. We consider the Roe algebra \cite{Ro} associated to $\Gamma$.  This  is the crossed product C$^\ast$-algebra $$\mathfrak R(\Gamma) = C^{\ast}(\Gamma \times l^{\infty}(\Gamma)),$$ where $\Gamma$ acts on $l^{\infty}(\Gamma)$ by left translations.

The algebra $\mathfrak R(\Gamma)$ admits a canonical representation, which we will denote by $\pi_{\rm Koop}$, into $\B(l^2(\Gamma))$. Generally, in the literature, it the  the C$^\ast$-algebra $\pi_{\rm Koop}(\mathfrak R(\Gamma))$ that  is referred to as to  the Roe algebra.   It is well known (see e.g. \cite{BO}) that the image C$^\ast$-algebra  is  the reduced $C^{\ast}$-algebra $C^{\ast}_{\rm red}(\Gamma \rtimes \ell^{\infty}(\Gamma))$.  If $\Gamma$ is exact the above  C$^\ast$- crossed product algebras are all isomorphic.

The partial action of $G$ on $\Gamma$ induces an action of $G$ by partial automorphisms on $\mathfrak R(\Gamma)$. We denote the corresponding groupoid crossed product $C^{\ast}$-algebra by
$$
\mathfrak R_G(\Gamma) = C^{\ast}(G \rtimes \mathfrak R(\Gamma)) =
C^{\ast}((G \rtimes \Gamma)\rtimes\ell^{\infty}(\Gamma)).
$$
We    consider the partial,  semidirect  product $G \rtimes \Gamma$ and let it act naturally, as  a groupoid  on $\ell^{\infty}(\Gamma)$.
If  $G = \Gamma$, we have $G \rtimes \Gamma= \Gamma  \times \Gamma^{\rm op}$, where $\Gamma^{\rm op}$ is
the same group as $\Gamma$, but with opposite multiplication. The group $\Gamma^{\rm op}$ acts by right multiplication on $\Gamma$ and hence on $\ell^{\infty}(\Gamma)$.  In this case, $\mathfrak R_G(\Gamma)$    is the C$^{\ast}$-algebra 
$$\mathfrak R_\Gamma(\Gamma)=C^{\ast}((\Gamma \times \Gamma^{\rm op}) \rtimes l^{\infty}(\Gamma)).$$ If $G$ is larger then $\Gamma$,  then  the C$^\ast$-algebra $\mathfrak R_G(\Gamma)$ contains
the crossed product C$^{\ast}$-algebra  $$C^{\ast}((G \times G^{\rm op}) \rtimes C(K)),$$ where $K$ is the profinite completion of $\Gamma$, with respect to the subgroups that are domains for the partial transformations in $G$ on $\Gamma$.

We denote by  $\pi_{\rm Koop} : \mathfrak R_G(\Gamma) \to \B(l^2(\Gamma))$ the canonical C$^\ast$-algebra representation extending the representation of  $\mathfrak R(\Gamma)$ into $\B(l^2(\Gamma))$. Consider the Calkin projection $\pi_{\rm Calk}$ of $\B(l^2(\Gamma))$ into the Calkin algebra $\mathcal Q(l^2(\Gamma))$.

 In this paper we analyse the  representation
 \begin{equation}\label{calkin}
\Pi_{\mathcal Q}=\pi_{\rm Calk} \circ \pi_{\rm Koop} :\mathfrak R_G(\Gamma) \to \mathcal Q(l^2(\Gamma)).
\end{equation}
We  find sufficient conditions so that $\Pi_{\mathcal Q}$ factorizes to a representation of  the reduced crossed product algebra
 $$\mathfrak R_{G, {\rm red}}(\Gamma) = C^{\ast}_{\rm red}(G \rtimes C^{\ast}_{\rm red}(\Gamma \times l^{\infty}(\Gamma)))=C^{\ast}_{\rm red}((G \rtimes \Gamma)\rtimes\ell^{\infty}(\Gamma)).$$
 

 This is equivalent to the fact that the $C^\ast$ algebra generated by left and right convolutors, by  elements in the group  $\Gamma$, acting on $\ell^2(\Gamma)$,  is   isomorphic, modulo the ideal of compact operators (\cite{Do}), to the reduced C$^\ast$-algebra
$C_{\rm red}^\ast(\Gamma\times\Gamma^{\rm op})$ associated to 
 the group $\Gamma\times\Gamma^{\rm op}$.

This property of the group $\Gamma$ 
is designated in the literature  (\cite{Oz}, \cite{CAD}) as the property AO.
This property was introduced by Akemann and Ostrand in \cite{AO}, where
they proved that the above property holds true for the free~groups. As noted explicitly  in (\cite{CAD}), the property AO is very close to the property $\mathcal S$ of N. Ozawa (\cite{Oz}). 

The properties AO and $\mathcal S$ were proven to hold true for much larger class of discrete
groups, first by Skandalis \cite{Sk}, where it was proven to hold true for lattices of Lie
groups of rank 1, and then by Ozawa \cite{Oz} for hyperbolic groups (see also \cite{HG}).

In this paper we take the point of view that the    states $\varphi$ on the image in the Calkin algebra of the crossed product $C^{\ast}$-algebra 
$
\mathfrak R_G(\Gamma)$ are  realised by computing the displacement, under the action of $G\rtimes \Gamma$, on  finite measure sets in infinite Loeb-measures (\cite{Lo}) spaces. 

We recall that in non-standard analysis (\cite{Ro}), given a countable   set $X$  and a free ultrafilter $\omega$ one constructs the associated non standard universe
denoted by $^\ast X$, consisting of all sequences in $X$ that coincide eventually, relative to the ultrafilter $\omega$. Let $\beta(X)$ be the Stone- \v Cech compactification of $X$. Then  $\beta(X)$ consists of all ultrafilters limits on $X$, that is $\beta(X)$ is the space of characters of $\ell^\infty (X)$. Then $X^\ast$  admits a canonical projection onto $\beta(X)$, associating to each sequence in $^\ast X$, the corresponding ultrafilter that sequence defines on $X$ (see e.g. \cite{Stroyan}, \cite{Li}).

 We consider the  subalgebra of subsets of $X$ consisting of internal subsets of $X^\ast$.
The internal subsets of $^\ast X$  are subsets obtained as follows. Let $A=(A_n)_{n\in \N}$ be a sequence of subsets of $X$. We denote by $\mathcal C_{\omega}(A)$ the internal subset defined by the formula
$$\mathcal C_{\omega}(A)=\mathcal C_{\omega}((A_n)_{n\in \N})=\{(x_n)_{n\in \N}\in  X^\ast\mid x_n \in A_n, {\rm eventually\  relative \ to\  } \omega \}.$$
The internal sets are closed with respect to finite reunions, intersections and set differences. In the case the sets $((A_n)_{n\in \N})$ are eventually finite (with respect to $\omega$), the set $\mathcal C_{\omega}(A)$ is called a hyperfinite set. Its cardinality is  $${\rm card\ } A=({\rm card\ } A_n)_{n\in \N})\in ^\ast \N.$$

Similar to the internal sets, the hyperfine sets form an algebra $\mathfrak A_{\rm hfinite}$. To  every sequence of finite support probabilities $\mu_n$ on $X$ one associates naturally a finitely additive measure on $\mathfrak A_{\rm hfinite}$, defined by the formula  
\begin{equation}\label{genloeb}
\mu_\omega=\mathop{\lim}_{n \to \omega} \mu_n.
\end{equation}

A breakthrough construction, due to P. Loeb (\cite{Lo}, (for a concise introduction see \cite{Li}, \cite{Cut},  \cite{El}, \cite{Pa}), proves that 
that the measure $\mu_\omega$ extends to the $\sigma$-algebra $\mathfrak A_{\rm L}$ generated by the hyperfinite sets. This can be seen as a consequence of the classical theorem of Caratheodory.   A more  constructive is obtained by using    the $\chi_1$  saturation principle (see e.g. \cite{Cut},\cite{Li}) that asserts that the reunion of a family $\mathcal H$  of hyperfinite sets is again a hyperfinite set if and only if the family $\mathcal H$ is finite.

Loeb construction shows that the  $\sigma$-algebra $\mathfrak A_{\rm L}$ to which the measure  $\mu_\omega$ extends to a countably additive measure, consists of subsets $Y\subseteq ^\ast X$ that may be "sandwiched", with arbitrary precision, with respect to the Loeb measure,  between sets in $\mathfrak A_{\rm hfinite}$. We refer to \cite{Li} for an excellent introductory  exposition of this construction. 
 
A remarkable example of Loeb measures as above, is obtained when  the measures $\mu_n$ have equal weights (equal to $\frac{1}{{\rm card (support(}\nu_n)}, n \in \N.)$ The corresponding Loeb measure is called a  Loeb uniform counting measure. It associated to  the cardinal
 $\alpha=  ({\rm card (support(}\nu_n))_{n\in \N})\in ^\ast \N.$ We denote this measure by $\nu_\alpha$ and note that it has an obvious extension to the $\sigma-$ algebra generated by hyperfinite sets of cardinality comparable to $\a$. Moreover $\nu_\a$ is obviously invariant to transformations on $^\ast X$ induced by permutations of $X$. We will prove that all Loeb measures as in formula (\ref{genloeb}) are absolutely continuous with respect to a direct sum of, mutually singular, Loeb uniform counting measures as above.  

Because of the work of Calkin (\cite{Ca}), it can be shown that 
states $\varphi$ of the type  considered above, obtained by restriction to $\mathfrak R_G(\Gamma)$ of states  defined  on the Calkin algebra, are limits of convex combinations of states of the form
$$\mathop{\lim}_{n\to \omega}\langle \cdot\ \mu_n^{1/2}, \mu_n^{1/2}\rangle,$$
with $(\mu_n)_{n\in \N}$ as in formula (\ref{genloeb}). To such a    state one associates      a canonical Loeb measure $\mu_\omega$. We use measures as above     to construct an explicit description of states $\varphi$ as above. A related  measure construction is the notion of a density measure (see \cite{Ma}, \cite{SZ}).

In general the Loeb measure associated with an arbitrary  state $\varphi$  is not a Loeb uniform counting measure, but it is  a limit of convex combinations of states corresponding to measures that are absolutely continuous with Loeb uniform counting measure of various corresponding to hyperfinite sets of various hyperfinite cardinality. 
A similar correspondence between   states and  Loeb measures  was introduced in   (\cite{El}, \cite{El1}). We refer to the above mentioned paper for an  extensive review of Loeb measure techniques.

 Every hyperfinite set $\mathcal C_\omega(A)$ as above determines a (generally infinite)  $G\rtimes \Gamma$-invariant measure space $\Y_A$. The corresponding  Loeb uniform counting measure $\nu=\nu_{{\rm card\ }A}$ determines a $G\rtimes \Gamma$-invariant measure (Lemma \ref {$y_a$}).  As in (Section 2.6 in \cite {El}, see  also \cite{El1},  \cite {mah1}) , we may assume that the $\sigma-$ algebra of measurable sets is separable, by restricting to the $\sigma$-algebra generated by the translates, by $G\rtimes \Gamma$, of $\mathcal C_\omega(A)$, that is also closed to the module action on $\Y_A$ of $\ell^\infty(\Gamma)$.  
 
  The crossed product algebra  $
\mathfrak R_G(\Gamma)$ has a natural embedding into the crossed product C$^\ast$-algebra $$C^{\ast}((G \rtimes \Gamma) \rtimes L^{\infty}(\Y_A, \nu_{{\rm card\ }A})).$$

Using this, we prove that the GNS representation associated to a state $\varphi$ as above is weakly contained in the direct sum of the Koopman representations (\cite{Ke}) of the form $$\pi_{{\rm Koop},A}: 
C^{\ast}((G \rtimes \Gamma) \rtimes L^{\infty}(\Y_A, \nu_{{\rm card\ }A})) \to
\mathcal B(L^{2}(\Y_A, \nu_{{\rm card\ }A})),$$
\noindent restricted to $
\mathfrak R_G(\Gamma)$.  We deduce the behaviour of states $\varphi$ as considered above from the analysis   of
the Koopmann representations of the enveloping C$^\ast$ algebras, introduced above, corresponding to hyperfinite sets $A$ that avoid any given finite subset of $\Gamma$.
Let $\beta(\Gamma)$ be the Stone-\v Cech compactification of $\Gamma$ and let 
$\partial(\beta(\Gamma))=\beta(\Gamma)\setminus \Gamma$ be its boundary.

The spaces $\Y_A$ are $G \rtimes \Gamma$ invariant subsets of the non standard universe  $^\ast \Gamma$. Hence there is a $G \rtimes \Gamma$-equivariant projection $\pi_{\Y_A}: \Y_A \to  \partial(\beta(\Gamma))$. 
If $K$ is a profinite completion of $\Gamma$ the above projection may be further composed with the projection onto $K$. We get a canonical projection $\pi_{K,A}:\Y_A\to K.$

We assume that $\Gamma$ is exact. It follows that the C$^\ast$-algebra $$
\mathfrak R_A(\Gamma)=C^{\ast}( \Gamma \rtimes L^{\infty}(\Y_A, \nu_{{\rm card\ }A}))$$ is nuclear. We analyse the action of $G$ on the center $\cZ$ of the von Neumann algebra  obtained by taking the weak closure  of the image, in the Koopman representation, of  $\mathfrak R_A(\Gamma)$.

Because of the nuclearity condition,  the analysis is reduced to the case when the action of $\Gamma$ (or of  a quotient of $\Gamma$ by an amenable subgroup) has a fundamental domain.
If the action of $G$ on $\Y_A$ is free then this allows to define a canonical $G$-invariant measure $\tilde \nu$ on the spectrum $\tilde \Y$ of $\cZ$.

 We prove that  there exists a choice of a hyperfinite set $\tilde A$, of the same hyperfinite  cardinality $\a$  as $A$, such that $G$-equivariantly $(\cZ,\tilde \nu)$
is equivalent to  $(\Y_{\tilde A}, \nu_{{\rm card\ }\tilde A})$. Moreover $\tilde A$ sits in the   fiber,  with respect to $\pi_{K,A}$, of the identity element in $K$. Here $G$ acts by conjugation. If all elements in $\Gamma$ determining elements in in the fiber have amenable stabilizers, it follows that the action of $G$ on $\cZ$ is amenable (\cite {CAD}, \cite{CAD1}, \cite{Oz1}). 

Consequently, if $G=\Gamma$,  the image through the Koopman representation $\pi_{{\rm Koop},A}$  of 
$C^{\ast}((G \rtimes \Gamma) \rtimes L^{\infty}(\Y_A, \nu_{{\rm card\ }A}))$ is nuclear (\cite{CAD}, \cite{CAD1}).  This forces the restriction of $\pi_{{\rm Koop},A}$ to the image of $
\mathfrak R_G(\Gamma)$ to factorize to  the reduced crossed product algebra.

The case when the action of $G$ on $\Y_A$ is non-free is treated separately, by easier arguments, related to the dynamics of the action of $G$ on cosets of amenable stabilizer groups.
For every $x\in \Gamma$ we denote by $$\{x\}'\cap \Gamma=C_\Gamma(x)=\{y\in \Gamma\mid yxy^{-1}= x\}$$ the centralizer of $x$ in $\Gamma$. 

Let $(\partial [\beta(\Gamma)])^{\rm fixed}$ be the subset of the boundary $\partial [\beta(\Gamma)]$
 of the Stone-\v Cech compactification
    obtained by taking   the closure, in $\beta(\Gamma)$, of the  reunion of
     all (necessary infinite) cosets of the form $C_\Gamma(x) y$,  $x\in \Gamma\setminus\{e\}, y\in \Gamma$. 
     
     Denote by $\mathcal A_0$ the set of conjugacy classes  of amenable subgroups of $\Gamma$.
     Using the above method we prove:

\begin{thm}\label{AOconditions}
 Let $\Gamma$ be a countable, exact (\cite{Ki},\cite {KW})  group with infinite non-trivial conjugacy orbits (i.c.c group). We assume that the set $\mathcal A_0$   is at most countable and that 
\item {(i)}There exists a family of finite index normal subgroups $(\Gamma_n)_{n\in \N}$, with trivial intersection, such that the conjugacy orbits of elements in $\Gamma$, with non amenable stabilizers, are separated in the profinite topology, defined by the family $(\Gamma_n)_{n\in \N}$, from the identity element.
\item{(ii)}
Assume that    $\Gamma \times \Gamma^{\rm op}$ acts amenably on $(\partial[\beta(\Gamma)])^{\rm fixed}$.

Then the  representation $\Pi_{\mathcal Q}$ of $C^{\ast}((\Gamma \rtimes \Gamma^{\rm op}) \rtimes \ell^{\infty}(\Gamma))$ into the Calkin algebra factorizes to a representation of  the reduced $C^{\ast}$ - algebra 
$$
C^{\ast}_{\rm red}((\Gamma \rtimes \Gamma^{\rm op}) \rtimes \ell^{\infty}(\Gamma)).
$$
\noindent In particular  it follows that $\Gamma$ has the ${\rm AO}$ property.
\end{thm}

As a corollary we obtain that both $\PGL_2(\Z[\frac{1}{p}])$ and $\SL_3(\Z)$ have the Akemann-Ostrand property (\cite{AO}). Recall that by \cite{Sa} these groups do not have the Ozawa's property $\mathcal S$ (\cite{Oz}, \cite{CAD}).
\begin{cor}
The groups $\PGL_2(\Z[\frac{1}{p}])$ and $\SL_3(\Z)$ have the AO property (\cite{Oz},\cite{CAD}).  Hence, using (\cite{Oz}), it follows that the group von Neumann algebras,
$\L(\SL_3(\Z))$ and $\L(\SL_n(\Z))$, $n\geq 4$, are non-isomorphic. 
\end{cor}

\section{Definitions and outline of the proof of main results}

We recall (\cite{Ca})   that a faithful family of  C$^\ast$-algebra representations  of the Calkin algebra  is obtained as follows. Let  $\omega$ be  a free ultrafilter on $\N$.  Let   $\xi=(\xi_n)_{n\in \N} $ be a sequence of norm 1 vectors (weights) in $l^2(\Gamma)$ that are  weakly convergent to 0. For each $n\in N$, let $\omega_{\xi_n, \xi_n}$ be the vector state $\langle \cdot\xi_n, \xi_n \rangle$ on $\B(l^2(\Gamma))$. Consider the weak limit state on $\mathfrak R_G(\Gamma)$
\begin{equation}\label{stateomega}
\varphi_{\xi} = \mathop{\lim}\limits_{n \to \omega} \omega_{\xi_n, \xi_n}|_{\mathfrak R_G(\Gamma)}.
\end{equation}
 Consequently  the direct sum of the GNS-representations (see e.g. \cite{Sak}) associated with the  all states   introduced above, is a faithful representation of the algebra $\mathfrak R_G(\Gamma)$.
By obvious density and linearity arguments,  is sufficient to assume that all the vectors  $\xi_n, n\in \N$, appearing in formula (\ref{stateomega})  have positive entries, and have     finite support in $\Gamma$.

To obtain a better picture of the GNS C$^\ast$-algebra representation associated with the states $\varphi_{\xi}$,
we  use  the non-standard universe  $^{\ast}\Gamma$ associated with the ultrafilter $\omega$ (\cite{Ro},\cite{Lo}). 

Any sequence of positive weights $\xi=(\xi_n)_{n\in \N}$ of   positive  weights on $\Gamma$, of finite support in $\Gamma$,  determines  a positive 
Loeb measure $\mu_\xi$ on the $\sigma$-algebra of Loeb measurable sets, constructed  as follows:

 Let $A=(A_n)_{n\in \N}$ be a sequence of finite sets in $\Gamma$. We denote by $\mathcal C_\omega(A)$  the corresponding hyperfinite subset of $^{\ast}\Gamma$.
Then we define
\begin{equation}\label{Loebw}
\mu_\xi(\mathcal C_\omega(A))=\mathop{\lim}\limits_{n \to \omega}\sum_{a \in A_n} \xi^n(a)
\end{equation}
Consider two weight sequences  $\xi_i=(\xi^n_i)_{n\in \N}, i=1,2 $, as above. Their pairing induces a finite Loeb measure on  $^{\ast}\Gamma$. The measure is defined by the requirement that   to a hyperfinite subset $A= (A_n)_{n\in \N}$, it associates  the value
$$
( \d \mu_{\xi_1}, \d \mu_{\xi_2} )(A )= \mathop{\lim}\limits_{n \to \omega} \sum_{a \in A_n} \xi_1^n(a)\xi^n_2(a).
$$
We denote the total mass of this measure by $\langle \d \mu_{\xi_1}, \d \mu_{\xi_2} \rangle$.
Clearly  this is computed as 
\begin{equation}\label{pairing}
\langle \d \mu_{\xi_1}, \d \mu_{\xi_2} \rangle = \mathop{\lim}\limits_{n \to \omega} \sum_{a \in \Gamma} \xi_1^n(a)\xi^n_2(a).
\end{equation}

Consequently, any sequence  $\xi=(\xi_n)_{n\in \N}$ of   positive  weights, of finite support,  corresponding to norm 1 vectors in $\ell^2(\Gamma)$,    determines a hyperfinite (\cite{Lo}) Loeb probability measure on $^{\ast}\Gamma$, defined by the formula
\begin{equation}\label{loebprob}
\nu_\xi= (\d \mu_{\xi}, \d \mu_{\xi} ).
\end{equation}
 The  support of this measure  is contained in the hyperfinite set $$\mathcal C_\omega(A)=\mathcal C_\omega( (A_n)_ {n\in \N})= \mathcal C_\omega(({\rm supp\ }\xi_n)_ {n\in \N}).$$
%
%

We endow the   space  $^{\ast}\Gamma$,  with the   $\sigma$ - algebra of Loeb  measurable sets with respect to hyperfinite all cardinalities.  Then $^{\ast}\Gamma$  is a  fibration over the Stone-\v Cech  compactification $\beta\Gamma$ of $\Gamma$ (\cite{Stroyan}).
Denote by $\B(^{\ast}\Gamma)$ the algebra of  bounded Loeb measurable functions on $^{\ast}\Gamma$. We have a  canonical embedding
\begin{equation}\label{ce}
\Phi : l^{\infty}(\Gamma)=C(\beta\Gamma)\to \B(^{\ast}\Gamma).
\end{equation}

Let $\Y_A$ be the minimal  subset of $^{\ast}\Gamma$,  closed to multiplication by functions in $\Phi(\ell^{\infty}(\Gamma))$ and to  translations by elements in  $G$ and $\Gamma$, containing the hyperfinite set $\mathcal C_\omega( (A_n)_ {n\in \N})$ corresponding to the sets $(A_n)_ {n\in \N}$. Let $\B_{}(\Y_A)$ be the algebra consisting of restrictions to $\Y_A$ of functions in  $\B(^{\ast}\Gamma)$. Then 
$\B_{}(\Y_A)$ has a canonical $l^{\infty}(\Gamma)$ module structure.
We denote by 
\begin{equation}\label{phi} 
\Phi_A : l^{\infty}(\Gamma) \to \B_{}(\Y_A),
\end{equation} 
the  representation obtained by restriction to $\Y_A$ of the embedding introduced in formula (\ref{ce}).
Recall that $\partial(\beta\Gamma)=\beta(\Gamma)\setminus \Gamma$.
If the sets $(A_n)_{n\in \N}$ avoid eventually, with respect to the ultrafilter $\omega$, any given finite subset of $\Gamma$, we obtain a representation, also
 denoted by $\Phi_A$, 
\begin{equation}\label{phir}
\Phi_A : C(\partial(\beta\Gamma)) \to \B_{}(\Y_A).
\end{equation}
In particular $\B_{}(\Y_A)$ is a bimodule over $C(\partial(\beta\Gamma))$.

We use the notation $C_0(\cdot)$ to denote the dense subalgebra of the crossed product  $C^{\ast}$ - algebra consisting of sums of finite support with respect the group elements.
Clearly $G \rtimes \Gamma$ acts by translations on $\B(^{\ast}\Gamma)$, and hence we get a algebraic representation of $C_0((G \rtimes \Gamma) \rtimes l^{\infty}(\Gamma))$ into $\B(^{\ast}\Gamma)$, which we denote by $\pi_{\rm Koop}$. 

 Fix $(g,\gamma)\in G\rtimes \Gamma$. Let $\chi_S$ be a characteristic function in $\ell^\infty(\Gamma)$. Denote by $
 (g, \gamma)(\d \mu_{\xi})$ the pullback of the measure $\d \mu_{\xi}$ by the transformation on $\Y_A$ induced by the partial transformation 
$(g,\gamma)\in G\rtimes \Gamma$. We consider  $\chi_S(g, \gamma)$ as an element in $C^{\ast}((G \rtimes \Gamma)\rtimes l^{\infty}(\Gamma))$. 

Then the value of the
 state $\varphi_{\xi}$, introduced in formula (\ref{stateomega}) at this given  element  in $C^{\ast}((G \rtimes \Gamma)\rtimes l^{\infty}(\Gamma))$, using formula (\ref{pairing})  is computed as follows.
\begin{equation}
\varphi_{\xi}[\chi_S(g, \gamma)] = \langle (g, \gamma)(\d \mu_{\xi}), \Phi(\chi_S)\d \mu_{\xi} \rangle.
\end{equation}

The  pairing of measures introduced in formula (\ref{pairing}) is not a Hilbert space scalar product. To obtain a scalar product,   one should find common measure, with respect to which the   measures $\langle (g,\gamma)(\d \mu_{\xi}), (g,\gamma)(\d \mu_{\xi})\rangle $, and hence the measures $\langle \d \mu_{\xi}, (g,\gamma)(\d \mu_{\xi})\rangle $, for all  $(g,\gamma)\in G\rtimes \Gamma$,   are simultaneously absolutely continuous.


Let $\a=(\alpha_n)_{n\in \N}\in\ ^{\ast}\N$ be  a non-standard natural number. Let $\nu_{\a}$ be  the Loeb measure on $^{\ast}\Gamma$, which compares the cardinality of hyperfinite subsets to $\a$ (see formula (\ref{uniform})).
Assume that  $A=(A_n)_{n\in \N}$ are finite subsets of $\Gamma$ such that
${\rm card\ } A=({\rm card\ } A_n)_{n\in \N}= \a$.

 Consider the  measurable subset $\Y_A\subseteq ^\ast \Gamma$ introduced above. Then $\nu_\a$    
restricts to  a $G\rtimes\Gamma$ invariant measure on $\Y_A$. One  constructs (see Lemma \ref{(phia)}) a canonical    representation 
$\Phi_ {A,\a}$, extending $\Phi_A$,  of the crossed product $\mathfrak R_G(\Gamma)$ into the crossed product C$^\ast$-algebra 
$C^{\ast}((G \rtimes \Gamma) \rtimes L^{\infty}(\Y_A, \nu_{{\rm card\ }A}))$.

 We introduce the following equivalence relation on $^\ast\N$. If $\a_1, \a_2 \in \ ^{\ast}\N$, then $$\alpha_1\cong \alpha_2 \rm{\  if  }  \mathop{\lim}\limits_{n \to \omega}\dfrac{\a_1(n)}{\a_2(n)}\in (0,  \infty).$$
 Clearly, if $\alpha_1, \alpha_2$ are non-equivalent  then $\nu_{\a_1}$ is singular to $\nu_{\a_2}$.   
 
  We choose a family of mutually singular measures $(\nu_{\a})_{\a \in \ ^{\ast}\N}$, corresponding to a choice of representatives in $\mathcal N\subseteq\ ^{\ast}\N/\cong$.  Then $\mathop{\oplus}\limits_{\a\in \mathcal N }\nu_{\a}$ is a $G\rtimes \Gamma$ - invariant measure on $^{\ast}\Gamma$. 
  
   We prove that for any probability measure $\nu_{\xi}$,  as above, there exists a countable subset $\mathcal N_0\subseteq \mathcal N $, such that the  direct sum of the countable family of, mutually singular measures  $\mathop{\oplus}\limits_{\a\in \mathcal N_0 }\nu_{\a}$ dominates $ \nu_{\xi}$, and hence also dominates $(g, \gamma)(\nu_{\xi})$ for all $(g,\gamma)\in G\rtimes \Gamma$.

The measures in the direct sum $\mathop{\oplus}\limits_{\a\in \mathcal N_0 }\nu_{\a}$ are mutually singular. Hence,  to determine the continuity properties of the representation of $\mathfrak{R}_G(\Gamma)$ into the Calkin algebra,  it  is sufficient to determine, with $(\Y_A, \nu_{{\rm card\ }A})$ as above,  the continuity properties of the  Koopmann representations $\pi_{\rm Koop,A}$  (see Definition \ref{koopmann}) of the enveloping $C^{\ast}$ - algebras
 $$C^{\ast}((G \rtimes \Gamma) \rtimes L^{\infty}(\Y_A, \nu_{{\rm card\ }A}))= C^{\ast}(G \rtimes C^{\ast}(\Gamma \rtimes L^{\infty}(\Y_A, \nu_{{\rm card\ }A}))),$$  into $\B(L^2(\Y_A, \nu_{{\rm card\ }A}))$.

We assume that $\Gamma$ is exact and that $G \rtimes \Gamma$ acts freely on $\Y_A$. The non-free case is treated separately in Section \ref{nonfree}. Since $^\ast \Gamma $ is a a measurable fibration over $\beta(\Gamma)$, it follows that  $\Gamma$ acts amenably on $\Y_A$. By this statement we understand  the fact that the commutative C$^\ast$-algebra $C^{\ast}(\Gamma \times L^{\infty}(\Y_A, \nu_{{\rm card\ }A}))$ is nuclear. The argument that we use here, following \cite{CAD1}, is the fact that the commutative  C$^\ast$-algebra $L^{\infty}(\Y_A, \nu_{{\rm card\ }A})$ is a central $\Gamma$-$C(\partial(\beta \Gamma))$ algebra and that the group  $\Gamma$ acts amenably (\cite{CAD1}, Proposition 8.2 and \cite{Oz1}) on 
$\partial(\beta \Gamma).$

  Assume that $\Gamma$ is a group with infinite, non-trivial conjugacy classes (briefly i.c.c.). Then, the center of the von Neumann algebra generated by the image of the crossed product $C^{\ast}(\Gamma \rtimes L^{\infty}(\Y_A, \nu_A))$, through  the Koopmann representation into $\B(L^2(\Y_A, \nu_{{\rm card\ }A}))$,  consists of the $\Gamma$ - invariant functions in $L^{\infty}(\Y_A, \nu_{{\rm card\ }A})$ (\cite{strat}).

The center coincides with the center $\mathcal Z$ of the von Neumann algebra generated by the reduced $C^{\ast}$ - algebra $C^{\ast}_{\rm red}(\Gamma \rtimes L^{\infty}(\Y, \nu_{{\rm card\ }A}))$. The later algebra is constructed, through GNS construction from the semifinite trace  induced by   the $\Gamma$ - invariant measure $\nu_{{\rm card\ }A}$.
We analyze the center $\mathcal Z$ in the context of the $C^{\ast}_{\rm red}$ - algebra representation. The corresponding von Neumann algebra, generated by $C^{\ast}_{\rm red}(\Gamma \rtimes L^{\infty}(\Y, \nu_{{\rm card\ }A}))$,  is either of type I or of type II.

We prove that the first case corresponds to the existence of a fundamental domain for the action of $\Gamma$ on $\Y_A$. In the second case, there exists an amenable subgroup $\Gamma_0$ in $\Gamma$, that invariates a measurable subset $F_0$ of $\Y_A$, such that the action of $\Gamma$ on $\Y_A$ is obtained by Mackey's induction (\cite{mackey}) from the action of $\Gamma_0$ on $F_0$.
In this case, the translates of the subset $F_0$, by distinct elements on $\Gamma / \Gamma_0$, are mutually disjoint (modulo zero measure sets).
In either case we transfer the measure from $F$ or $F_0$ onto $\mathcal Z$. We  obtain that  $\mathcal Z$ is of the form $L^{\infty}(\tilde{\Y}, \tilde{\nu})$, where $\tilde{\nu}$ is a $G$ - invariant measure on $\tilde{\Y}$.

 We assume that we are given a family of normal, finite index subgroups $\Gamma_n$ of $\Gamma$, that separate points in $\Gamma$. In the case $G$ is different from $\Gamma$, the subgroups in the family are associated to the subgroup lattice associated with the subgroups that appear as domains of the partial isomorphisms defining the action of $G$ on $\Gamma$. Let $K$ be the profinite completion of $\Gamma$ with respect to this family of subgroups. We assume that non-trivial  elements  in $\Gamma$ having non-amenable stabilizers under the conjugacy action of $G$ are separated in the profinite topology from the identity element in $\Gamma$.
 
   Consider the embedding $C(K) \subseteq l^{\infty}(\Gamma)$. We compose this with the   representation  $\Phi_A$ of $l^{\infty}(\Gamma)$ into $L^{\infty}(\Y_A, \nu_{{\rm card\ }A})$  introduced in formula (\ref{phi}). We obtain a representation $\Phi_{K,A}: C(K) \to L^{\infty}(\Y_A, \nu_{{\rm card\ }A})$. Denote the corresponding  projection $\pi_{K,A} : \Y_A \to K$.
   In the case of type II, $K$ will be replaced by a quotient.
   
    Using the measure $\tilde \nu $ introduced above we prove that  the center  $\mathcal Z$ may be realised  as a subalgebra of  the bounded measurable  functions defined on  the fiber  $\pi_{K,A}^{-1}(\{e\})$.
In this identification we prove that     the measure $\tilde{\nu}$, and the corresponding action of $G$ on $\tilde{\Y}$, may be modelled as $(\Y_{\tilde A}, \nu_{{\rm card\ }\tilde A})$   using a different hyperfinite subset $(\tilde{A}_n)_n$ of $^{\ast}\Gamma$.

The subsets  in the family   $\tilde A=(\tilde{A}_n)_{n\in N}$ are obtained by translating by elements in $\Gamma$ pieces of the  sets  $A_n$. The requirement is that the sets   $\tilde{A}_n$ are contained in the subgroups $\Gamma_{k_n}$, $n\in \N$, for 
a suitable choice of the sequence $(k_n)_{n\in \N}$ of integers, tending to infinity, so that the corresponding translations are behaving similarly to the corresponding translations of the  sets in the ultrafilter limit after $\omega$.

The assumption that the orbits of elements in $\Gamma$, that have non-amenable stabilizer subgroups under the conjugation action of $G$, are separated from the identity in the profinite topology, implies that the action of $G$ on $(\Y_{\tilde{A}}, \nu_{\a})$ is amenable (\cite{CAD1}). Hence so is the action of $G$ on $\mathcal Z$. Because $\mathcal Z$ is the center of $$\{\pi_{\rm Koop,A}(C^{\ast}(\Gamma \rtimes L^{\infty}(\Y_A))\}'' ,$$ it follows that the cross product C$^\ast$-algebra  $$C^{\ast}(G \rtimes \pi_{\rm Koop,A}(C^{\ast}(\Gamma \times L^{\infty}(\Y_A, \nu_A)))$$ is nuclear.

 It follows that crossed product C$^\ast$-algebra  $$\pi_{\rm Koop,A}(G \rtimes C^{\ast}(\Gamma \rtimes L^{\infty}(\Y_A, \nu_A)))$$ is nuclear, and hence $\pi_{\rm Koop,A}\circ \Phi_{A,\a}$
  is a representation of the reduced $C^{\ast}$ - algebra $C^{\ast}_{\rm red}(G \rtimes C^{\ast}(\Gamma \times l^{\infty}(\Gamma))$.


We analyze separately the case of points in $^{\ast}\Gamma$ having non-trivial stabilizers under the action of $G \rtimes \Gamma$. It is easily seen that such points are contained in the image (under the embedding of $l^{\infty}(\Gamma) \subseteq L^{\infty}(\Y_A)$) of the subset $(\partial(\beta\Gamma))^{\rm fixed}$  of $\partial(\beta\Gamma)$ obtained as a reunion of the cosets of stabilizer  subgroups under the conjugation action.
In the case of $G = \Gamma = \PGL_2(\Z[\frac{1}{p}])$ the stabilizer groups are cyclic.

As the intersection of maximal abelian subgroups of $\Gamma$ is trivial, the cosets of maximal abelian subgroups of $\Gamma$ have finite intersection, and hence their image in $\partial(\beta\Gamma)$ have trivial intersection. Thus on $(\partial(\beta\Gamma))^{\rm fixed}$ the action of $G \rtimes \Gamma = \Gamma \times \Gamma^{\rm op}$ is weakly contained in $l^2(\Gamma \times \Gamma^{\rm op} /\Gamma_0 \times x\Gamma_0 x^{-1})$ where $x \in \Gamma$ and $\Gamma_0$ is one of the stabilizers.
A similar argument works for $\SL_3(\Z)$.

\section{Loeb measures and states on $\mathfrak{R}_G(\Gamma)$}

We analyze the structure of Loeb measures ([Lo], [Li]) with weights.
We prove that every hyperfinite probability Loeb measure is absolutely continuous with respect to the direct sum of uniform Loeb measures, corresponding to various cardinalities.

We use this to construct $G\rtimes \Gamma $-infinite measure spaces $(\Y, \nu)$, whose associated Koopman representations  are in turn  used to represent the states of the $C^{\ast}$ - algebra $\mathfrak{R}_G(\Gamma)$. To prove that the states obtained this way exhaust the states on $\mathfrak{R}_G(\Gamma)$ corresponding to  the representation into the Calkin algebra $\cQ(l^2(\Gamma))$ representation, we use Calkin faithful representation of the algebra $\cQ(l^2(\Gamma))$,  determined by the choice, of an arbitrary free ultrafilter $\omega$.

States on $\cQ(l^2(\Gamma))$, when restricted to $l^{\infty}(\Gamma)\subseteq \B(^{\ast}\Gamma)$ define  Loeb measures, with weights, on the non-standard universe $^{\ast}\Gamma$.
To obtain a $\Gamma$ - invariant setting, we prove that such measures   are absolutely continuous with respect to uniform Loeb counting measures.

Denote by $\B(^\ast\Gamma)$, the algebra of functions, measurable with respect to the Loeb - $\sigma$ algebra $\mathfrak A_L$ generated by hyperfinite sets. We are implicitly proving that the element in $H^1(\Gamma, \B(^\ast\Gamma))$ corresponding to a Loeb weighted measure is trivial.

Let $\omega$ be a free ultrafilter on $\Gamma$ (which in the next theorem is considered just as a copy of $\N$).
To state the first result we will forget about the group structure on $\Gamma$. 
If  $(A_n)_{n\in \N}$ is a  sequence of  finite subsets of $\N$ we denote by
 $\cC_{\omega}((A_n)_{n \in \N})$ be the corresponding hyperfinite set in the non-standard universe $^{\ast}\N$ associated to the free ultrafilter $\omega$.
 We assume that $(A_n)_{n\in \N}$ avoids eventually any finite subset of $\N$.

Let $\xi=(\xi^n)_{n\in \N} $ be a sequence of probability measures on $\N$, with finite support $B_n$.
Let $\nu_\xi=\nu_{\xi, \omega}$ be the internal Loeb measure corresponding to $\xi$.
Thus 
$$\nu_{\xi, \omega} (\cC_{\omega}((A_n)_{n \in \N})=\mathop{\lim}\limits_{n \to \omega}\sum_{a \in A_n\cap B_n} \xi^n(a).$$
For $\alpha \in ^{\ast}\N$, we consider the Loeb measure on the  $\sigma$ - algebra $\B_{L}$, generated by hyperfinite sets, given by the formula
\begin{equation}\label{uniform}
\nu_{\alpha}(\cC_{\omega}((A_n)_{n \in \N})) = \mathop{\lim}\limits_{n \to \omega} \dfrac{{\card} \; A_n}{\a_n}.
\end{equation}

Clearly $\nu_{\alpha}$ is supported on the $\sigma$-algebra of  hyperfinite sets of (non-standard) hyperfinite cardinality equal to $\a$. We let $\cN_0$ be $^{\ast}\N$, modulo the equivalence relation
$$
\a \sim \beta \; {\rm if} \; \mathop{\lim}\limits_{n \to \omega} \dfrac{\a_n}{\beta_n} \in (0, \infty).
$$
\noindent Clearly $\nu_{\a}$, $\nu_{\beta}$ are mutually singular if $\a \neq \beta$ in $\cN_0$.
We have

\begin{thm}\label{disintegration}
Given $\nu_{\xi, \omega}$ as above, there exists internal functions $(f_{\a})_{\a \in \cN_0}$ on $^{\ast}\N$ (only a countable number different from 0) such that for every $\epsilon >0$, there exists a measurable subset $A_\epsilon$ of measure $\nu_{\xi, \omega}(A_\epsilon)>1-\epsilon$, such that
$$
\nu_{\xi, \omega}|_{A_\epsilon} =
 \mathop{\oplus}\limits_{\a \in \cN_0} f_{\a}\d\nu_{\a}|_{A_\epsilon}.
$$
\end{thm}

\begin{proof}
We proceed by maximality.
Assume that we found a chain, $(f_{\a}, \nu_{\a})_{\a \in \cN}$, so that the internal positive functions  $f_{\a}$  have  support in $\cC_{\omega}((A_n^{\a})_{n \in \N})$, where
$$\alpha_n={\rm card \ }(A_n^\a), n \in \N, \a \in \cN_0,$$  and such that
$$
\nu_{\xi, \omega}\mid_{\cC_{\omega}((A_n^{\a}))} = f_{\a}\d\mu_{\a}, \a \in \cN_0.
$$

Assume $\nu_{\xi, \omega}(\cC_{\omega}((A_n)_{n \in \N}) \setminus (\mathop{\bigcup}\limits_{\a \in \cN_0} \cC_{\omega}((A_n^{\a})_{n \in \N})) \neq 0$.

In the complement, by $\chi_1$ - exhaustation principle we find $\cC_{\omega}((A_n^0)_{n \in \N})$, such that $\nu_{\xi, \omega}\mid_{\cC_{\omega}((A_n^0)_{n \in \N})} \neq 0$, and
$$
{\rm supp}\nu_{\xi, \omega}\mid_{\cC_{\omega}((A_n^0)_{n \in \N})} = \cC_{\omega}((A_n^0)_{n \in \N}).
$$

Let $A_n^{0, M} = \{ a \in A_n^0 \mid \xi_n(a) \leq \dfrac{M}{\card A_n}\}$.

Then, by the support condition
$$
\nu_{\xi, \omega}(\mathop{\cup}\limits_{M} \cC_{\omega}((A_n^{0, M})_{n \in \N})) \neq 0
$$

\noindent and hence at least one of the sets  $\cC_{\omega}((A_n^{0, M_0})_{n \in \N})$ has non-zero $\mu_{\xi}$ measure.
Then $\nu_{\xi, \omega} \mid_{\cC_{\omega}((A_n^{0, M_0})_{n \in \N})}$ is absolutely continuous with respect $\mu_{(\card A_n^0)_n}$.

Hence, by maximality and  transfinite induction we obtain an infinite chain as above.
 Since the measure $\nu_{\xi, \omega}$ evaluated at the sets $\cC_{\omega}((A_n^{\a})_{n \in \N})$ is non-zero, the process will exhaust the support of $\nu_{\xi, \omega}$ after a countable number of steps.
 
 This completes the proof, except for the fact that elements in $\cN_0$ may eventually be repeated in the maximal chain. Again using the $\chi_1$-saturation principle, we may replace a reunion of disjoint  elements of the form \break ${\cC_{\omega}((A_n^{\a,s})_{n \in \N})}$, $s \in S$ corresponding to a unique class $\a \in \cN_0$, by a single
 set of the same form, with an arbitrary small loss in measure with respect to the measure $\nu_{\xi}$. 
 Since the set of possible $\alpha$ is countable, using the above process, and taking each time an approximation of the order $\epsilon/2^n, n \in \N$ we obtain the result.
\end{proof}

In the case when the countable set $\N$ is replaced by $\Gamma$, we also have the action of the  $G\rtimes \Gamma$, by left translations, on $^{\ast}\Gamma$.
In this case given $\a = (\a_n) \in \ ^{\ast}\N$, and given $\cC_{\omega}(A_n)$ a hyperfinite subset of $^{\ast}\Gamma$ of cardinality $\a = (\a_n)$, we introduce the following definition:
\begin{defn}\label{$y_a$}
The $G\rtimes \Gamma$ - invariant measure space associated to $(\cC_{\omega}(A_n))$ is the measure space $(\Y_A, \nu_{\a})$, where
$\alpha_n={\rm card \ } (A_n), n\in \N,$ and
$$
\Y_A =\Phi_A(\ell^\infty(\Gamma)) \mathop{\bigcup}\limits_{(g, \gamma)\in G\rtimes \Gamma } \cC_{\omega}[((g,\gamma){A_n})_{n \in \N}] \subseteq \ ^{\ast}\Gamma.
$$
\noindent This space is endowed with the $G\rtimes \Gamma$ - invariant Loeb measure $\nu_{\a}$ associated to $\a \in\ ^{\ast}\N$.

Unless $(A_n)_{n \in \N}$ is contained in sequence of Folner sets for the group $\Gamma$, which is precluded if $\Gamma$ is non-amenable, it will follow that $(\Y_A, \nu_{\a})$ is an infinite measure space. Remark that the measures $\nu_{\a}, \nu_{\beta}$ remain mutually singular if $\a \neq \beta$ in $\cN_0$.
\end{defn}

\begin{lemma}\label{(phia)} We use the notations introduced above.
Assume $\Gamma$ is contained as an almost normal subgroup on a larger group $G$. Then $G$ acts by partial isomorphisms on $\Gamma$ and hence on $^{\ast}\Gamma$, and on $\B(^{\ast}\Gamma)$. Consequently we have a natural $C^{\ast}$ representation
$$
\Phi_{A, \omega}: C^{\ast}(G \rtimes C^{\ast}(\Gamma \rtimes l^{\infty}(\Gamma))) \rightarrow C^{\ast}(G \rtimes C^{\ast}(\Gamma \rtimes L^{\infty}(\Y_A, \nu_{\a}))).
$$
\end{lemma}

\begin{proof}
The above representation is constructed as follows.
First note that $G$ acts on $^{\ast}\Gamma$, and this is compatible with the action of $\Gamma$.
Let $\B(^{\ast}\Gamma)$ be the space of bounded functions on $^{\ast}\Gamma$ that are measurable with respect to the Loeb $\sigma$ - algebra generated by hyperfinite sets.

We use the representation $\Phi_A:\ l^{\infty}(\Gamma) \to  L^{\infty}(\Y_A, \nu_{\a})$  introduced in formula (\ref{phi}). Since  the above inclusion is $\Gamma$ and $G$ equivariant, it passes to the crossed product algebra.

We recall that we are  using  the notation $C_0(\cdot)$ to denote the dense subalgebra of the crossed product  $C^{\ast}$ - algebra consisting of sums of finite support with respect the group elements.

At the algebraic level we have
$$
C_0^{\ast}(G \rtimes C_0^{\ast}(\Gamma \rtimes l^{\infty}(\Gamma))) \subseteq C_0^{\ast}(G \rtimes C_0^{\ast}(\Gamma \rtimes \B(^{\ast}\Gamma))).
$$

Since $\Y_A \subseteq \B(^{\ast}\Gamma)$ is invariant under $G$ and $\Gamma$ it follows that the algebraic inclusion introduced  above also  gives a $C^{\ast}$ - representation
$$
C^{\ast}(G \rtimes C^{\ast}(\Gamma \rtimes l^{\infty}(\Gamma))) \rightarrow C^{\ast}(G \rtimes C^{\ast}(\Gamma \rtimes L^{\infty}(\Y_A, \nu_{\a}))).
$$
\end{proof}

\section{Representation of the algebra $C^{\ast}(G \rtimes C^{\ast}(\Gamma \rtimes l^{\infty}(\Gamma)))$}

In this section we use the embedding from Lemma \ref{(phia)} to determine the representations of the algebra $\mathfrak{R}_G(\Gamma)$ which appear in   the Calkin algebra  representation considered in formula (\ref{calkin}).

First we introduce  the following extended  definition  of the Koopmann representation.

\vspace{0.3cm}

\begin{defn}\label{koopmann}
 Given a group $H$ acting by measure preserving transformations, on a space $(X, \mu)$, let the Koopmann $C^{\ast}$ - representation (see e.g. \cite{Ke})
$$
\pi_{\rm Koop}: C^{\ast}(H \rtimes L^{\infty}(X, \mu)) \rightarrow \B(L^2(X, \mu))
$$

\noindent be the representation obtained by letting $H$ act by left translation on $L^{\infty}(X, \mu)$ and hence on $L^2(X, \mu)$.
Let $L^{\infty}(X, \mu)$ act canonically by multiplication on $L^2(X, \mu)$.

It is well known that  the two representations of $H$ and of the algebra $L^{\infty}(X, \mu)$ (\cite{Pe}) induce a representation of the  crossed product algebra $C^{\ast}(H \times L^{\infty}(X, \mu))$. 
\end{defn}

\

If $H$ doesn't preserve the measure, but only preserves the class of $\mu$, then using the cocycle $\T(g, x) = \dfrac{g(\d\mu)}{\d\mu}(x)$, $x \in X$, $g \in G$ one can still define a unitary representation of $H$ on $L^2(X, \mu)$ (using the cocycle $\T(g, x)$ to perturb the formula of $\pi_{\rm Koop}\mid_H$. It is known (\cite{Krieg})) that this also  extends to a unitary representation of $C^{\ast}(H \times L^{\infty}(X, \mu))$.

\begin{ex}
We endow $l^{\infty}(\Gamma)$ with the counting measure
$$
\varepsilon = {\rm Tr}_{\B(l^2(\Gamma))}\mid_{l^{\infty}(\Gamma)}.
$$

\noindent Then the embedding $\pi_{\rm Koop}$ $C^{\ast}(\Gamma \times l^{\infty}(\Gamma)) \subseteq \B(l^2(\Gamma))$ is exactly the Koopmann embedding.
This obviously extends to a representation, also denoted by $\pi_{\rm Koop} $  of the crossed product algebra $C^{\ast}(G \rtimes C^{\ast}(\Gamma \times l^{\infty}(\Gamma)))$, since the action of  $G$ invariates $l^{\infty}(\Gamma)$.
\end{ex}

Let $\pi_Q : \B(l^2(\Gamma)) \to \cQ(l^2(\Gamma))$ be the Calkin representation.
Then the states determining the topology on  $\cQ(l^2(\Gamma))$ are obtained as follows:

\begin{lemma}\label{Cal} (\cite{Ca})
Let $\omega$ be a free ultrafilter on $\N$. Let $\xi_n \in l^2(\Gamma)$ be a sequence of vectors converging weakly to 0. By linearity and density we may assume that the vectors $\xi_n$ have finite  support in $\Gamma$, with positive entries  and that $(\xi_n)^2$ they are probability measures on $\Gamma$.

Then  the weak limit  state
\begin{equation}\label{state}
\varphi_{\xi,\omega} = \mathop{\lim}\limits_{n \to \omega} \langle \cdot \xi_n, \xi_n \rangle
\end{equation}
factorizes to a state on $\cQ(l^2(\Gamma))$.
The states $\varphi_{\xi,\omega}$  introduced above, through the Gns-Representation determine the topology on $\cQ(l^2(\Gamma))$.
\end{lemma}

\begin{defn}
Consider  two Loeb weighted measures $\mu_{\xi, \omega}$, $\mu_{\eta, \omega}$ as in formula (\ref{Loebw}). We define $(\mu_{\xi, \omega}, \mu_{\eta, \omega})$ as the measure $\dfrac{\d\mu_{\eta, \omega}}{\d\mu_{\xi, \omega}}\mu_{\xi, \omega}$. This is  a measure on $^{\ast}\Gamma$ (when if not absolutely continuous we take the above to be zero).

Then $\langle \mu_{\xi, \omega}, \mu_{\eta, \omega} \rangle$ is by definition $\int(\mu_{\xi, \omega}, \mu_{\eta, \omega})\d\mu_{\xi, \omega}$.
This is $$\mathop{\lim}\limits_{n \to \omega} \mathop{\sum}\limits_{\gamma \in \Gamma} \xi^n(\gamma)\overline{\eta^n(\gamma)}.$$
\end{defn}

The following result was noted in the introductory section.

\begin{lemma}
Consider  the restriction of the  state $\varphi_{\xi, \omega}$ on $\B(l^2(\Gamma))$, introduced in formula (\ref{state}) to $\pi_{\rm Koop}(C^{\ast}(G \rtimes (C^{\ast}(\Gamma \rtimes l^{\infty}(\Gamma))))) \subseteq \B(l^2(\Gamma))$.

Then for $(g, \gamma) \in G \rtimes \Gamma$, $\chi_S \in l^{\infty}(\Gamma)$ we have,  using the notations from the  previous definition:
\begin{equation}\label{eqphi}
\varphi_{\xi, \omega}((g, \gamma)\chi_S) = \langle\Phi(\chi_S)\mu_{\xi, \omega}, g^{-1}(\mu_{\xi, \omega})\rangle
\end{equation}

\end{lemma}

We use the decomposition from Theorem \ref{disintegration} from the previous chapter.

\begin{thm}

Let  $A = (A_n)_{n \in \N} $  be a sequence of finite subsets of $\Gamma$ that eventually avoids, with respect to the ultrafilter $\omega$, any  finite subset of $\Gamma$.  Let 
\begin{equation}\label{alpha}
\alpha=(\alpha_n)_{n \in \N}=({\rm card\ } A_n)_{n \in \N}.
\end{equation}
 Consider the  the representation $\Phi_{A, \omega}$ constructed in Definition \ref{(phia)}. 
Let $(\Y_A, \nu_{\a})$ be a measure   space  as introduced in Definition \ref{$y_a$}. Let $\pi_{\rm Koop,A}$ be the associated  Koopmann representation (see Definition \ref{koopmann}):
 $$\pi_{\rm Koop,A}: C^{\ast}(G \rtimes C^{\ast}(\Gamma \rtimes L^{\infty}(\Y_A, \nu_{\a})))\to \B(L^2(\Y_A, \nu_{\a})).$$

Then the direct sum of the representations of $\mathfrak{R}_G(\Gamma)$  of the form  $$ \pi_{{\rm Koop}, A} \circ \Phi_{A, \omega} $$  weakly contains  the representation $\Pi_{\mathcal Q}=\pi_{\rm Calk} \circ \pi_{\rm Koop}$, introduced in formula (\ref{calkin}), of    $\mathfrak{R}_G(\Gamma)$ in $\cQ(l^2(\Gamma))$.
\end{thm}

\begin{proof}

Consider a state on $\mathfrak{R}_G(\Gamma)$ as in  Lemma \ref{Cal}. By Theorem \ref{disintegration},
the $\rm GNS$ representation associated with the state $\varphi_{\xi, \omega}$ on $\mathfrak{R}_G(\Gamma)$, is weakly contained in the direct sum of the $C^{\ast}$ - representation of $\mathfrak{R}_G(\Gamma)$ obtained as follows:
We decompose  $\nu_{\xi, \omega}$  as the direct sum of the  mutually singular measures
$$
\oplus f_{\a}\nu_{\a}\mid_{\cC_{\omega}((A_n^{\a}))}
$$

\noindent and using the formula (\ref{eqphi}), we get that the state $\varphi_{\xi, \omega}$ is a limit of convex combinations of states of the form
$$
x \in R_G(\Gamma) \to \langle \pi_{{\rm Koop}, A}(\Phi_{A,\omega}(x)) f^{1/2}_{\a}, f^{1/2}_{\a} \rangle_{L^2(\Y_{A}, \nu_{\a})}
$$
\end{proof}

\begin{cor} To prove that the representation $\pi_{\rm Calk} \circ \pi_{\rm Koop, A}$ of $\mathfrak R_G(\Gamma)$ factorizes to the reduced crossed product $\mathfrak R_{G,{\rm red}}(\Gamma)$,
it is  sufficient to prove that  the  states of the form
$$
\varphi_{\omega, A}(x)= 
  \langle \pi_{{\rm Koop}, A}(\Phi_{A,\omega}(x)) \chi_F, \chi_F \rangle_{L^2(\Y_A, \nu_{{\rm card} A})}, 
$$
\noindent where $\chi_F$ is the characteristic function corresponding to the hyperfinite set $\cC_{\omega}(A_n)$ determining the space $(\Y_A, \nu_{{\rm card} A})$, are continuos states on the reduced C$^\ast$-algebra $\mathfrak R_{G,{\rm red}}(\Gamma)$.
\end{cor}
\begin{proof}

This is proved by taking the Koopmann representation $\pi_{{\rm Koop}, A}$ of $$C^{\ast}(G \rtimes C^{\ast}(\Gamma \times L^{\infty}(\Y_A, \nu_{\a})))$$ and restrict it to the image of $\mathfrak{R}_G(\Gamma)$ through the embedding $\Phi_{A, \omega}$.

\end{proof}

In the next chapter we analyze the algebra $C^{\ast}(G \rtimes (C^{\ast}(\Gamma \rtimes L^{\infty}(\Y_A, \nu_{\a}))))$.
Because  $\nu_{\a}$ is a $(G \rtimes \Gamma)$ invariant measure, it follows that the reduced $C^{\ast}$-algebra $C^{\ast}_{\rm red}(G \rtimes (C^{\ast}(\Gamma \rtimes L^{\infty}(\Y_A, \nu_{\a}))))$,  has weak closure  which is a   type I or II von Neumann algebra.  On the other hand the von Neumann algebra  $$\{\pi_{{\rm Koop}, A} \circ \Phi_{A, \omega}(\mathfrak{R}_G(\Gamma))\}''\subset  \B(L^2(\Y_A, \nu_{\a})),$$ could be a-priori of any von Neumann type (\cite{Sak}).


\begin{rem}
The fact that the cocycle in $H^1(\Gamma, \Y_A)$ corresponding to the measure defined up to equivalence measure from $\mu_{\xi, \omega}$ on $\Y_A$ vanishes simplifies the proof.
Otherwise one should consider on $\Y_A$ also measures that are $(G \rtimes \Gamma)$ - equivariant, up to equivalence of measures.
\end{rem}

\section{Structure of the center of the von Neumann algebra $\{ \pi_{{\rm Koop},A}(C^{\ast}(\Gamma \rtimes L^{\infty}(\Y_A, \nu_{\a}))) \}''$}

In this section we assume that $\Gamma$ is exact.  Consider a measure space $(\Y_A, \nu_{\a})$, as  constructed in the previous section. 
We prove below the $C^{\ast}$ - algebra $C^{\ast}(\Gamma \rtimes L^{\infty}(\mathcal Y_A, \nu_{\a}))$ is nuclear.  We prove that the center splits into a type I or type II component, corresponding to the behaviour of the action of $\Gamma$ on $\Y_A$. We also prove that  the center $$\cZ=\cZ\big(\{ \pi_{\rm Koop}(C^{\ast}(\Gamma \rtimes L^{\infty}(y_A, \nu_{\a}))) \}''\big)$$ carries a canonical measure $\tilde{\nu}$. In the next chapter we prove that  the measure $\tilde{\nu}$ is invariant under the action of $G$.

We assume in the following sections, until Section \ref{nonfree},  that $\Gamma$ acts freely on $\Y_A$. The non-free case will be treated separately in Section \ref{nonfree}. 
To distinguish between the two cases,  we let $(^\ast \Gamma)^{\rm fixed} $ be the subset of $^\ast \Gamma$, consisting of points in $^\ast \Gamma$ having non-trivial stabilizers, relative to the action of $G\rtimes \Gamma$.
Let
\begin{equation}\label{free}
(^\ast \Gamma)^{\rm free}=^\ast \Gamma\setminus (^\ast \Gamma)^{\rm fixed}.
\end{equation}

\begin{lemma}\label{nuclear} Let $(\Y_A, \nu_{\a})$ be as in the previous section. 
Then
 the crossed product  C$^\ast$- algebra  $C^{\ast}(\Gamma \rtimes L^{\infty}(\Y_A, \nu_{\a}))$ is nuclear.
\end{lemma}

\begin{proof}
Indeed, as we noted in the Section \ref{intro} , we have a $\Gamma$ - equivariant (in fact $G \rtimes \Gamma$ equivariant), measurable projection
\begin{equation}\label{cech}
\pi_{^{\ast}\Gamma}\  :\  ^{\ast}\Gamma \to \beta\Gamma.
\end{equation}

Restricting $\pi_{^{\ast}\Gamma}$  to $\Y_A$, we obtain a measurable projection 
\begin{equation}\label{pia}
\pi_{\Y_A} : \Y_A \to \partial(\beta\Gamma),
\end{equation}
 that is $\Gamma$ - equivariant and measurable (see \cite{Stroyan}, \cite{Li}).
 
 Then $L^{\infty}(\Y_A, \nu_{\a}))$ is a central $\Gamma$-$C(\partial(\beta(\Gamma)))$ algebra and by exactness $\Gamma$ acts amenably on $\partial(\beta(\Gamma))$.

It follows that (\cite{CAD}) that 
 the crossed product $C^{\ast}$ - algebra $C^{\ast}(\Gamma \rtimes L^{\infty}(\Y_A, \nu_{\a}))$ is nuclear.

\end{proof}

As a consequence of Lemma \ref {nuclear} we deduce that $$\pi_{\rm Koop,A} : C^{\ast}(\Gamma \rtimes L^{\infty}(\Y_A, \nu_{\a}))\to \mathcal B(L^2(\Y_A, \nu_{\a}))$$ factorizes to $C^{\ast}_{\rm red}(\Gamma \rtimes L^{\infty}(\Y_A, \nu_A))$. The later $C^{\ast}$-algebra is easily represented, using he fact that $\nu_A$ is invariant under $\Gamma$ and hence defines a semifinite trace on $C^{\ast}_{\rm red}(\Gamma \rtimes L^{\infty}(\Y_A, \nu_A))$,  also denoted by $\nu_A$. This semifinite trace  is used to construct the GNS-representation denoted by $\pi_{{\rm red, A}}$. The Hilbert space associated to the GNS representation is $L^2(\L(\Gamma \rtimes L^{\infty}(\Y_A, \nu_A))), \nu_A)$.

\begin{lemma}We use the notations introduced above.  Assume that the measure space $\Y_A$ is contained in $ (^\ast \Gamma)^{\rm free}$.
We also  assume that $\Gamma$ is an i.c.c. group.  Let  $\cZ$ be the center of
the von Neumann algebra $$\{ \pi_{{\rm Koop},A}(C^{\ast}(\Gamma \rtimes L^{\infty}(\Y_A, \nu_{\a}))) \}'' \subseteq \B(L^2(\Y_A, \nu_{\a})).$$
Then $\cZ$  is isomorphic to  the center of the von Neumann algebra generated by the $C^\ast_{\rm red}$ representation of  $C^{\ast}(\Gamma \rtimes L^{\infty}(\Y_A, \nu_A))$. It is identified with the $\Gamma$ - invariant functions in $L^{\infty}(\Y_A, \nu_{\a})$.
\end{lemma}

\begin{proof}
In the case of the Koopmann representation this follows simply because $L^{\infty}(\Y_A, \nu_{\a})$ is maximal abelian in $\B(L^2(\Y_A, \nu_{\a}))$. In the reduced $C^{\ast}$ - algebra case this follows from the i.c.c. condition (see e.g. \cite{strat}).
\end{proof}

To analyze the center we will  use the reduced $C^{\ast}$ - algebra representation and the identification of the center with the algebra of $\Gamma$-invariant, bounded measurable functions. This is algebra is determined by the analysis of the corresponding fundamental domains.
%
In the next definition we introduce a specific terminology to describe the analogue of a fundamental domain for the action of a discrete group, relative to a subgroup. 

\begin{defn}
Assume that $\Gamma$ acts freely by measure preserving transformations on a measure space $(\cX, \mu)$. Let $\Gamma_0$ be a subgroup of $\Gamma$. Let $F$ be a subset of $\cX$ that is invariant under $\Gamma_0$ and assume that if $s_1\Gamma_0 \neq s_2\Gamma_0$ then the intersection of $s_1F$ and $s_2F$ has zero measure.

If the above condition holds true we will say  that $F$ is a $\Gamma/\Gamma_0$ wandering subset of $\cX$. If, in addition, the translates of $F$ cover $\cX$, we will say that $F$ is a $\Gamma/\Gamma_0$ - wandering and generating subset of $\cX$ with respect to the action of $\Gamma$.

Note that in this case $\gamma\Gamma_0\gamma^{-1}$ invariates $\gamma F$, and $\gamma F$ is a $\gamma\Gamma_0\gamma^{-1}$ - wandering, generating subset. 
\end{defn}

\

\

The action of $\Gamma$ is (up to a choice of representatives) determined by action of $\Gamma_0$ on $F$. Indeed this follows from Mackey's construction (\cite{mackey}) that we recall below.

\begin{defn}[Induced action, \cite{mackey}] \label{mackey}%
We use the  context of the previous  definition.
The induced action on $(\Gamma \setminus \Gamma_0) \times F$  is constructed as follows:

\noindent Let  $s_i\Gamma_0, s_j\Gamma_0$, be two cosets of $\Gamma_0$ in $\Gamma$ and let $\gamma\in \Gamma$  such that 
 $$\gamma(s_i\Gamma_0) = s_j\Gamma_0.$$
This is equivalent to the fact there exists $\T(\gamma, s_i) \in \Gamma_0$ such that $$\gamma{s_i} = s_j \T(\gamma, s_i).$$
%
Let  $f \in F$.
Then, we define  $$\gamma(s_i\Gamma_0, f)=(s_j\Gamma_0, \T(\gamma, s_i)f) \in \Gamma \setminus \Gamma_0 \times F.$$

The fact that the above formula defines   an action of $\Gamma$ on $\Gamma/\Gamma_0 \times F$ follows directly from Mackey construction (\cite{mackey}) of the induced representation.
Then, the Koopmann unitary representation of $\Gamma_0$ on $L^2(F, \mu)$ induces the Koopmann unitary representation on $L^2(\cX, \mu)$, where $\cX=\Gamma \setminus \Gamma_0 \times F$ and the action of $\Gamma$ is as introduced above.
\end{defn}

\

We analyze below the decomposition of the centre $\cZ$,
\begin{equation}\nonumber
\cZ=\cZ\big(\big\{\pi_{{\rm red},A}( C^{\ast}(\Gamma \rtimes L^{\infty}(\Y_A, \nu_A))\big\}''\big)\subseteq \B(L^2(\L(\Gamma \rtimes L^{\infty}(\Y_A, \nu_A))), \nu_A)),
\end{equation}
 according to the type of action of $\Gamma$.

Recall that we denote by $\mathcal A_0$ the set of conjugacy classes  of amenable subgroups of $\Gamma$ and 
that we are  assuming that $\mathcal A_0$ is at most countable. 

\begin{thm}\label{center} Let $(\Y_A, \nu_{\a})$, $\cZ$ as above. There exists a decomposition of  $\cZ$ and $\Y_A$ subject to the conditions (i), (ii) below.

The center $\cZ$ is divided as a direct sum $$\cZ = \cZ_I \oplus \cZ_{II},\ 
\cZ_{II}= \mathop{\oplus}\limits_{\Gamma_0 \in \mathcal A_0} \cZ_{\Gamma_0}.$$ The decomposition of the center induces a    corresponding decomposition    $$\Y_A=\Y_I \oplus \Y_{II}, \ \Y_{II} = \mathop{\oplus}\limits_{\Gamma_0 \in \mathcal A_0} \Y_{\Gamma_0}.$$ Here $\Y_I, \Y_{II}$ are measurable subsets of $\Y_A$, that are  $\Gamma$ - invariant, with measure zero overlaps. Similarly $\Y_{\Gamma_0}, \Gamma_0\in \mathcal A_0$ are  $\Gamma$ - invariant, measurable subsets of $\Y_{II}$  forming a measurable partition of $\Y_{II}$.   

 The above decomposition has   the following properties.

\item {(i)}
There exists a measurable subset $F$ of $\Y_I$ that is a   $\Gamma$ - wandering, generating subset of $\Y_I$.

\item{(ii)} Let  $\Gamma_0\in \mathcal A_0$. 
Then there exists a measurable $\Gamma \setminus \Gamma_0$-wandering, generating subset $F_{\Gamma_0}$ of $\Y_{\Gamma_0}$. Up to a choice of representatives for cosets of $\Gamma \setminus \Gamma_0$, the action of $\Gamma$ on $\Y_{\Gamma_0}$ is $\Gamma$-equivariantly equivalent to the  action of $\Gamma$, obtained by Mackey's induction (Definition \ref{mackey}) from the action of $\Gamma_0$ on $F$.
\end{thm} 

\begin{proof}
The weight $\nu=\nu_{\a}$ is semifinite, and   $\Gamma$ acts by measure preserving transformations on $\Y=\Y_A$, which is a $G\rtimes \Gamma$-equvariant,  measurable fibration over  the non-discrete  spectrum  $\partial (\beta(\Gamma))=\beta(\Gamma)\setminus \Gamma$ of $l^{\infty}(\Gamma)$. As observed above   the algebra $C^{\ast}_{\rm red}(\Gamma \rtimes L^{\infty}(\Y, \nu))$ is nuclear because of the exactness of $\Gamma$.

We have a canonical semifinite trace on this algebra, obtained as the composition of the canonical, normal conditional expectation $E$ onto $L^{\infty}(\Y, \nu)$ with the infinite measure (weight) on $L^{\infty}(\Y, \nu)$ given by $\nu$. We consider the Koopman unitary representation $\pi_{\rm Koop}$ of the  crossed product $C^{\ast}$ - algebra $$C^{\ast}(\Gamma \rtimes L^{\infty}(\Y, \nu)),$$ on the Hilbert space $\H_{\nu}=L^2(\Y, \nu)$ associated to the semifinite trace $\nu$. We denote this C$^\ast$-algebra of  by $$C^{\ast}_{\rm Koop}(\Gamma \rtimes L^{\infty}(\Y, \nu))\subseteq B(L^2(\Y, \nu)).$$
\noindent  Because of nuclearity, the representation $\pi_{\rm Koop}$ of $C^{\ast}(\Gamma \rtimes L^{\infty}(\Y, \nu))$ into $C^{\ast}_{\rm Koop}(\Gamma \rtimes L^{\infty}(\Y, \nu))$    is isometric.

Let $$M=\overline{C^{\ast}_{\rm Koop}(\Gamma \rtimes L^{\infty}(\Y, \nu))}^{\rm w}\subseteq B(L^2(\Y, \nu)) ,$$ be the corresponding von Neumann algebra, which is necessary of semifinite type. Let $D = L^{\infty}(\Y, \nu)$ be the corresponding MASA in $M$, and let $E$ be the normal conditional expectation from $M$ onto $E$. Because of the infinite conjugacy classes condition on the group $\Gamma$, the center $\cZ=\mathcal{Z}(M)$ is is contained in $D = L^{\infty}(\Y, \nu)$. As observed above, $\cZ$ consists of the $\Gamma$-invariant functions in $D$.

 We identify the algebra $\mathcal{Z}(M)$  with  the algebra $L^\infty(\tilde{\Y}, \tilde{\nu})$, for some measure space $\tilde{\Y}$, and  a canonical measure $\tilde{\nu}$ on $\tilde{\Y}$  introduced  below. We denote the subsets of $\tilde \Y$ corresponding to summands $\cZ_I$, $\cZ_{II}$ and $\cZ_{\Gamma_0}$, $\Gamma_0\in \mathcal A_0$ respectively  by 
 $\tilde{\Y}_I$, $\tilde{\Y}_{II}$ and $\tilde{\Y}_{\Gamma_0}$
 
   In the case of type I  the measure $\tilde{\nu}$ is defined simply by letting $\tilde{\nu}(\tilde F)=\nu (F)$, if $F$ is a minimal  measurable subset of $\Y$, of finite measure such that the characteristic function  $\chi_{\tilde F}$ is the central support  in $M$ of the projection $\chi_F$.  In the case of type II, in the $\Y_{\Gamma_0}$ component, $\Gamma_0\in \mathcal A_0$, we impose the additional requirement that $F$ be a subset left invariant by $\Gamma_0$.
 
  The measure $\tilde{\nu}$ is in fact, in the case of type I , the Plancherel measure of the corresponding type I algebra (\cite{Dix}, Section 18). In the case of type II, the measure $\tilde{\nu}$ is the obvious analogue of the Plancherel measure. Moreover, as explained above, in both cases,  $L^\infty(\tilde{\Y}, \tilde{\nu})$ is the $\Gamma$-invariant part of $L^{\infty}(\Y, \nu)$.  
 
We denote by $\nu$, the semifinite, faithful weight on $M$ induced by $\nu_{\a}$. By nuclearity  $M$ can only have type $I_{\infty}$ or hyperfinite type $II_{\infty}$ components.  The fact that we get only type I$_\infty$ or  II$_\infty$ components  is a consequence of the absence of Folner sets. Indeed, by the nuclearity of the algebra $C^{\ast}(\Gamma \rtimes L^{\infty}(\Y, \nu))$, the type $II$ components are hyperfinite (\cite{Co}).

We disintegrate $M$ over the center $\mathcal{Z}(M)$. We obtain, almost everywhere with respect $\tilde{\nu}$, fibers $M_z \supseteq D_z$, $z\in \tilde{\Y}$, endowed a with normal faithful conditional expectation $E_z : M_z \to D_z$.  Since we have not given yet  the complete definition of the measure $\tilde \nu$, the notion of $\tilde\nu$-a.e. refers here to any measure $\tilde\nu$ which gives the isomorphism $\mathcal{Z}(M)\cong L^\infty(\tilde{\Y}, \tilde{\nu})$.

By disintegration over the center, the semifinite trace $\nu$, yields, for $z \in \tilde{\Y}$, almost everywhere, a semifinite trace  $\nu_{z}$  on $D_z$, extending to  a semifinite faithful trace on $M_z$.

In the case of type $I$, which corresponds to the central part $\Y_I$, because of the existence of a normal conditional expectation onto the algebra $D_z$, it follows that the algebras $D_z$ are maximal abelian, diagonal algebras. Hence any field of minimal projections is the multiplication operator on $L^{2}(\Y, \nu)$, with the  characteristic function a fundamental domain for the action of $\Gamma$

In the case of type $II$, which corresponds to the $\Y_{II}$ part in the statement, because of the fact that there exists a conditional expectation from $M_z$ onto $D_z$,  and since $M_z$ is of type $II_{\infty}$ it follows that $M_z$ admits a splitting $$M_z\cong N_z \otimes \B(l^2(I_z)),$$ where $N_z$ is a type $II$, (hyperfinite) factor, and $l^2(I_z)$ is the Hilbert space associated to a countable set $I_z$ (a.e. for $z\in \tilde\Y_{II}$). Below, we denote by  $\nu_0^z$  the canonical trace on $M_z$ obtained by the disintegration of $\nu$.

Moreover, since $D_z$ is maximal abelian and generated by finite projections, it follows that $D_z$ splits as $D_z^1 \otimes D_z^2$, in such a way that $D_z^1$ is a MASA in $N_z$ and $D_z^2$ is the maximal abelian diagonal algebra of $\B(l^2(I_z))$ associated to the basis indexed by $I_z$.

Let $\pi_z$ be the disintegration of the left regular representation of the group $\Gamma$ in $\H_{\nu}$.
Thus $$\{\pi_z(\Gamma), D_z\}'' = M_z, $$ for $z\in \tilde\Y_{II}$ almost everywhere. Then the unitary operators $\pi_z(\gamma)$ normalise the algebra $D_z$ for every $\gamma$. Consequently, there exists a permutation $P_z(\gamma)$ of $I_z$, $P_z(\gamma) : I_z \to I_z$ such that if $(e_{i, j}^z)$ is the matrix unit of $\B(l^2(I_z))$ associated to the basis indexed by $I_z$, then there exists unitary operators $u_i^z(\gamma)$, $i \in I_z$ in the normaliser $\mathcal N_{N_z}(D_z^1)$, such that $$\pi_z(\gamma) = \mathop{\sum}\limits_{i \in I_z} u_i^z(\gamma)\otimes e^z_{i, P_z(\gamma)(i)},\ \gamma \in \Gamma.$$
 Then, necessary, the map $\gamma \to P_z(\gamma)$ into the permutation group of $I_z$ is a homeomorphism. Hence there exists a subgroup $(\Gamma_0^z)$ of $\Gamma$ such that the index set $I_z$ is identified with the set of cosets $$\mathcal C_z =\{[s\Gamma_0^z]| s \in \Gamma\},$$ in $\Gamma / \Gamma_0^z$, $s\in \Gamma$, (a.e. for $z\in \mathcal Z$).   The above matrix unit is therefore  indexed by 
$\Gamma / \Gamma_0^z$.  
 We use the following notation for the matrix unit:   $$(e^z_{[s^z_i\Gamma^z_0],[s^z_j\Gamma^z_0}])_{[s^z_i\Gamma^z_0], [s^z_j\Gamma^z_0]\in \mathcal C_z} ,$$
 for $z$ almost everywhere.
 
  The above identification of the index set $I_z$ is $\Gamma$-invariant. The permutation $P_z(\gamma)$, in this identification, is translation by $\Gamma$ on $\Gamma/ \Gamma_0^z$, $\gamma \in \Gamma$. 
 Note that $\Gamma_0^z$ is necessary infinite, since otherwise we are back in the  case of type $I_\infty$.   
 
  Let $e_0^z$ in $\B(l^2(\Gamma/ \Gamma_0^z))$ be the projection corresponding to $e_{[\Gamma_0^z], [\Gamma_0^z]}$.
Then $e_0^z$ is fixed by $\pi^z(\gamma)$, $\gamma \in \Gamma_0^z$. Hence,  identifying $N_z$ with $N_z \otimes e_0^z$, we obtain a representation $\pi_0^z(\gamma)$, $\gamma \in \Gamma_0^z$ of $\Gamma_0^z$ into the unitary group of  $N_z$, such that the original representation $\pi^z$ is in this identification the induced representation ${\rm Ind}_{\Gamma_0^z}^{\Gamma}(\pi_0^z)$ on $$L^2(N_z, \nu_0^z) \otimes l^2(\Gamma/ \Gamma_0^z),$$ \noindent a.e. for $z\in \mathcal Z$.

Because in the original representation $E_z(\pi^z(\gamma)) = 0$, it follows that, if we denote by $\nu_0^z = \nu^z(e_0^z\cdot)$ the trace induced by $\nu$ on $N^z$, then $$\nu_0^z(\pi_0^z(\gamma)) =0, \ \gamma \in \Gamma_0^z\backslash\{e\}.$$ \noindent Moreover $\pi_0^z(\Gamma_0)'' \subseteq N_z$ and hence $N_z$  contains the type II$_1$ factor associated to the group $\Gamma_0^z$. The corresponding left regular  representation of $\Gamma_0^z$  normalises the Cartan subalgebra $D_z^1$, (a.e. for $z\in \mathcal Z$). Consequently, for $z$ almost everywhere, the factor  $N_z$ is the reduced crossed product von Neumann algebra
$$\mathcal L(\Gamma^z_0\rtimes D_z^1).$$

Since $N_z$ is hyperfinite, it follows that $\Gamma_0^z$ is amenable and infinite. Since $e_0^z$ is the projection in $D_z$ corresponding to $1\otimes e_{[\Gamma_0^z], [\Gamma_0^z]}$ it also follows  that the $\Gamma $-system $\Y_z$ (the fiber of $\mathcal Y$ at $z$ in the type II case) is $\Gamma$-equivariantly isomorphic to a $\Gamma$-system of the form $$\Y_z\cong F^z \times \Gamma / \Gamma_0^z,$$ where $F^z$ is a probability measure space, which is $\Gamma_0^z$ invariant and $\Gamma / \Gamma_0^z$ is endowed with the counting measure, almost everywhere for $z$ in $\mathcal Z$.  Since we have an at most  countable set of infinite amenable subgroups, the property (2) holds true, a.e. for $z\in \mathcal Z$.


\vspace{0.2cm}

 The subgroup  ${\Gamma}_0$ is uniquely determined, by construction, up to conjugacy by the algebra in the fiber. Hence  the sets $\Y_{\Gamma_0}$ are disjoint, when $\Gamma_0$ runs in $\mathcal A_0$ (up to overlaps of zero measure). 
\end{proof}

\vspace{0.2cm}

In the rest of the section we construct a Plancherel type  measure on $\cZ$ and give an explicit formula for its computation.

\vspace{0.2cm}

First we define a canonical measure $\tilde{\nu}$ on the center $\cZ$. The measure will be proven to be  canonical, in the sense that it is invariant under $G$, as we prove in the next section.

\begin{lemma}We use the previous notations and definitions.
The following construction defines a   measure $\tilde \nu$ on the center $\cZ$. Consequently, the center $\cZ$, will be  endowed with the  trace $\tilde{\nu}$. Denote  its spectrum by    $\tilde{\Y}$. This will be endowed with the measure $\tilde{\nu}$. 

\item{(i)}
 Because of the existence of a fundamental domain $F$ for the action of $\Gamma$ on $\Y_I$, the center $\cZ_I$, which consists   of  measurable, bounded $\Gamma$ - invariant functions,
 is canonically identified to bounded measurable functions on $F$.
 The measure $\nu_\alpha|_F$ induces the measure $\tilde \nu|_{\cZ_I}$ on 
$\cZ_I$. Denote by $\tilde{\Y}_I$ the corresponding subset of $\tilde{\Y}$.
\item{(ii)}
In the case of type II, let $\Gamma_0\in \mathcal A_0$. Let $F_{\Gamma_0}$ be a $\Gamma/\Gamma_0$ - wandering, generating subset of $\Y_A$, that is also $\Gamma_0$ invariant.
Then $\cZ_{\Gamma_0}$ is identified, as a measure space, with  $F_{\Gamma_0}$.  We transport the measure $\nu_\alpha|_{F_{\Gamma_0}}$
  onto $\cZ_{\Gamma_0}$.  The resulting measure is   $\tilde\nu|_{\cZ_{\Gamma_0}}$. The corresponding part of $\tilde {\Y}$ is denoted by $\tilde {\Y}_{\Gamma_0}$.

Consequently $\cZ=  L^{\infty}(\tilde{\Y}, \tilde{\nu})$, $\cZ_I=L^{\infty}(\tilde{\Y_I}, \tilde{\nu})$, 
$\cZ_{\Gamma_0}=
L^{\infty}(\tilde{\Y}_{\Gamma_0},
 \tilde{\nu}). 
$

\end{lemma}
\begin{proof} In the case of type $I$ this is obvious by construction. In the type II case corresponding to subgroup $\Gamma_0\in \mathcal  A_0$, we note that because we are using a disintegration process, it follows that in the fibers the von Neumann algebra 
 $\{C^{\ast}(\Gamma_0 \rtimes (F_0)_z)\}''$, for $z$ a.e in the spectrum of $  \cZ_{\Gamma_0}$, is a factor. Hence the $\Gamma$ invariant functions are again identified  with the  functions on a $\Gamma_0\setminus \Gamma$ wandering generating subset.

\end{proof}

We will use in the sequel a sequence of normal subgroups $\{ \Gamma_n \}_n$ of $\Gamma$, that have trivial intersection. In the case when $G \neq \Gamma$, we will use the groups in the family $\Gamma_g = g\Gamma g^{-1} \cap \Gamma$, $g \in G$, to construct the sequence.

We denote by $K$ the profinite completion of $\Gamma$ with respect to this family of subgroups.
Since $C(K) \subseteq l^{\infty}(\Gamma)$, we have that $\partial(\beta\Gamma)$ projects by a canonical, $G\rtimes \Gamma$-equivariant, continuous  projection $\pi_K$ onto $K$.

Let $\pi_{\Y_A}: \Y_A \to \partial(\beta\Gamma)$ be the canonical projection obtained by restriction from the canonical projection $\pi_{^{\ast}\Gamma} : \ ^{\ast}\Gamma \to \beta\Gamma$ (see the formulae (\ref{cech}),(\ref{pia})).
Then $$\pi_{K, A}=\pi_K \circ \pi_{\Y_A},$$
 is a $G\rtimes \Gamma$-equivariant, measurable projection from $\Y_A $ onto $ K$. We analyse the measurable structure of the fiber at $e$, with respect to $\pi_{K, A}$,   of the $\Gamma$ - invariant subsets of  $\Y_A$.

In the case of type II, we replace $\Gamma_n$ by $\tilde{\Gamma}_n = \Gamma_n\Gamma_0$. Then, the groups $\tilde{\Gamma}_n$ intersect exactly in $\Gamma_0$.
By the  normality of $\Gamma_n, n \in \N$, we obtain that
 $$ \tilde{\Gamma}_n = \Gamma_0\Gamma_n=\Gamma_n\Gamma_0.$$ 
\noindent  We let $K_{\Gamma_0}$ be the profinite limit of the coset spaces $\Gamma / \tilde{\Gamma}_n$, and in this case we obtain, similar to the above construction,   a projection $\pi_{K_{\Gamma_0}, A} : \Y_A \to K_{\Gamma_0}$. In the case of type II   we analyze the fiber at $[\Gamma_0] \in K_{\Gamma_0}$.

To do the analysis of the  of the trace of  $\Gamma$-invariant measurable sets in the fiber at $e$, we prove first an explicit formula for the intersection of sets in $\tilde{\Y}$. This is used to analyse the matrix coefficients of the action of $G$ on $\tilde{\Y}$.

\begin{defn}
Let $\Gamma_0$ be a subgroup as in part (ii) of Theorem \ref{center}. Let $s_i^n$ be  a sequence of coset representatives for $\Gamma_n$ in $\Gamma$ (respectively of $\tilde{\Gamma}_n$ in $\Gamma$).

Let $\overline{s_i^n\Gamma_n}$, (respectively $\overline{s_i^n\tilde{\Gamma}_n}$) be the profinite closure of the corresponding cosets in $K$ (respectively $K_{\Gamma_0}$). Let $C_i^n$ be $\pi^{-1}_{K, A}(\overline{s_i^n\Gamma_n}) \subseteq \Y_A$, (respectively $\pi^{-1}_{K_{\Gamma_0, A}}(\overline{s_i^n\tilde{\Gamma}_n})\subseteq \Y_A$). 
\end{defn}

\

We have the following formulae.

\begin{prop}\label{intersection}
We use the above notations. Let $\tilde{F}_0, \tilde{F}_1$ be measurable subsets of $\tilde{\Y}$, of finite measure.
Assume that $\chi_{F_0}$, $\chi_{F_1}$ are minimal projections (with the additional requirement that they are invariant to $\Gamma_0$ in the case of type $\cZ_{\Gamma_0}$, $\Gamma_0\in \mathcal A_0$) having central support $\chi_{\tilde{F}_0}$, $\chi_{\tilde{F}_1} \in \cZ$.
Clearly, in this case, $F_0, F_1$ are $\Gamma_0 / \Gamma_1$ - wandering subsets (simply $\Gamma$-wandering in the type I case).
Then
$$
\tilde{\nu}(\tilde{F}_0 \cap \tilde{F}_1) = \mathop{\sum}\limits_{\gamma \in \Gamma / \Gamma_0} \nu(F_0 \cap \gamma F_1) = 
$$
$$
= \mathop{\lim}\limits_{n \to \infty} \mathop{\sum}\limits_{i, j} ([(s_i^n)^{-1}[C_i^n \cap F_0]] \cap [(s_j^n)^{-1}[C_j^n \cap F_1]]).
$$
\noindent In the case of type I, we take $\Gamma_0 = \{e\}$ in the above formulae.
The same type of formula is valid for $n + 1$ sets $\tilde F_0, \tilde F_1, \ldots, \tilde F_n$. 
\end{prop}

\begin{proof}
The first part of the above equality follows from the $\Gamma / \Gamma_0$ - wandering property of $F_0, F_1$. The second property is an easy consequence from the first one, as when $n \to \infty$, the sets $$\{ (s_j^n)(s_i^n)^{-1} F_0 \cap F_1 \mid i,j = 1, \ldots, [\Gamma : \tilde{\Gamma}_n] \}$$ exhaust (biunivocally because of $\Gamma / \Gamma_0$-wandering property) the set of intersections $\gamma F_0 \cap F_1$.
\end{proof}

The last term in the equality in the above formula has the advantage that it may be treated as a formula for a measure in the fiber at $e$ of $\pi_{K, A}$, respectively $(\pi_{K_{\Gamma_0}, A})$.
The statement  will  be made more precisely when constructing the action of $G$ on $\cZ = L^{\infty}(\tilde{\Y}, \tilde{\nu})$.

\section{The action of $G$ on the center $\cZ = L^{\infty}(\tilde{\Y}, \tilde{\nu})$ and its generalised matrix coefficients}

In this section we construct a canonical action of the group $G$ on the center $\cZ=L^{\infty}(\tilde{\Y}, \tilde{\nu})$ and prove that this action is measure preserving. We use the notations and assumptions  from the previous sections. Assume  that $\Y_A\subseteq (^\ast \Gamma)^{\rm free}$.

We use the formula in Proposition \ref{intersection}  to describe the matrix coefficients of the  representation induced  by the Koopman  representation of \break $C^{\ast}(G \rtimes (C^{\ast}(\Gamma \rtimes L^{\infty}(\Y_A, \nu_{\a}))))$,  in the fiber of $\pi^{-1}_{K, \Y_A}(e)$.
This will be used to construct a different representation of the action of $G$ on $\cZ$, which under additional conditions will be amenable.

In the case $G = \Gamma$,  the group  $G$ acts by conjugation on $\Gamma$. Through the representation
$$
\pi_{\rm Koop,A} \circ \Phi_{A,\omega} : C^{\ast}(G \rtimes (C^{\ast}(\Gamma \times l^{\infty}(\Gamma)))) \rightarrow 
\B(L^2(\Y_A, \nu_{\a}))
,
$$

\noindent the group $G$ acts by the Koopmann representation, which we denote by $\pi^G_{{\rm Koop},A}$,  of $G$ on $L^2(\Y_A, \nu_{\a})$. Then $\pi^G_{\rm Koop}(G,A)$ normalises $L^{\infty}(\Y_A, \nu_{\a})$. 

Hence $\pi^G_{{\rm Koop},A}(G)$ normalises $$\pi_{{\rm Koop},A}(C^\ast(\Gamma \times L^{\infty}(\Y_A, \nu_{\a}))),$$ and hence normalises $\cZ$.

In general, if $\Gamma \subseteq G$, $\Gamma \neq G$, then $\Gamma$ is almost normal. Recall that $G$ acts (partially) by conjugation on $\Gamma$. The domain of the action of $g\in G$ is the subgroup $\Gamma_g$. This action is extended to an action of $G$ on $\Y_A$. We will use the conjugation notation to designate this action. Thus if $g\in G$ and $y\in \Y_A$, we  denote the action of $g$ on $y$ by $gyg^{-1}$.  We use the representation $\Phi_{A,\omega}$ introduced in Lemma \ref{(phia)}.  We have:

\begin{lemma}
Let $g \in G$, and let  $\chi_{\Gamma_{g^{-1}}},\chi_{\Gamma_g}\in \ell^{\infty}(\Gamma)$ be  the characteristic functions corresponding to the subgroups $\Gamma_g$, $\Gamma_{g^{-1}}$.  Then  $\pi_{\rm Koop}(g)$ is a partial isometry with initial space $\Phi_{A, \omega}(\chi_{\Gamma_{g^{-1}}})$ and range $\Phi_{A, \omega}(\chi_{\Gamma_g})$.

In the $C^{\ast}$ - algebra $C^{\ast}(G \rtimes (C^{\ast}(\Gamma \times l^{\infty}(\Gamma))))$, we have  the relation $$g(\gamma\chi_{\Gamma_{g^{-1}}})g^{-1} = (g\gamma g^{-1})\chi_{\Gamma_g}.$$
\noindent  This relation  is transferred, by  the   representation $\pi_{\rm Koop,A}$ into  $\B(L^2(\Y_A, \nu_A))$.
%
%
Hence  $g \in G$ will map $$\Phi_{A, \omega}(\chi_{\Gamma_{g^{-1}}})(\mathcal L(\Gamma_{g^{-1}} \rtimes L^{\infty}(\Y_A, \nu_A)))\Phi_{A, \omega}(\chi_{\Gamma_{g^{-1}}})$$
 \noindent  into $$\Phi_{A, \omega}(\chi_{\Gamma_g})(\mathcal L(\Gamma_g \rtimes L^{\infty}(\Y_A, \nu_A)))\Phi_{A, \omega}(\chi_{\Gamma_g}).$$
 \noindent

Then,  the  action of $G$ on $\cZ$ by partial isomorphisms  extends to an action of $G$ on $\tilde \Y$, constructed as follows. For $g\in G$ let  $s_i$ be the coset representatives for $\Gamma_g$ in $\Gamma$. We define, for $z \in \tilde \Y$,
$$
\T(g)(z) = \mathop{\sum}\limits_{i} \pi_{\rm Koop}(s_i) g(\Phi_{A, \omega}(\chi_{\Gamma_{g^{-1}}})z)g^{-1}.
$$
\end{lemma}

\begin{proof}
 Note that $\Phi_{A, \omega}(\chi_{\Gamma_g})$ is in fact the characteristic function of the subset in $\Y_A$ obtained as $\pi^{-1}_{K,A}(\overline{\Gamma}_g)$, where $\pi_{K,A} :\Y_A \to K$ is the canonical projection (and similarly in the type II case).
 
We look at $\Gamma$ - invariant functions, which are the elements of $\cZ$ as germs of $\Gamma$ - invariant functions. The action of $G$ on germs of such $\Gamma$ - invariant functions is clearly well defined, eventually by replacing the groups $\Gamma_g$ with normal subgroups $\Gamma_g^0 \subseteq \Gamma_g$, of finite index.
\end{proof}

We assume that the almost normal subgroup $\Gamma$ of $G$ verifies the additional condition that 
\begin{equation}\label{hecke}
[\Gamma : \Gamma_g] = [\Gamma : \Gamma_{g^{-1}}], g \in G.
\end{equation}

Using this condition in the case $\Gamma\ne G$, we obtain:
\begin{prop}
The above   defined action  of $G$ on $\cZ = L^{\infty}(\tilde{\Y}, \tilde{\nu})$ is measure preserving.
\end{prop}

\begin{proof}
Recall that $G$ acts by conjugation on $C^{\ast}(\Gamma \times l^{\infty}(\Gamma))$.

In the case $G = \Gamma$ this statement is obvious because: in the case of type I, since $G$ normalises $\Gamma$, it will map $\Gamma$ - wandering sets into $\Gamma$ - wandering set. Similarly, in the case of type II, it  maps $\Gamma / \Gamma_0$ wandering subsets into $\Gamma / g\Gamma_0 g^{-1}$ wandering subsets of $\Y_A$. Since the measure on $\cZ$ is obtained by transfer of the measure from the corresponding $\Gamma$ - wandering subsets it  follows that $G$ will preserve the measure on $\cZ$.

In the case when  $\Gamma$ is an almost normal subgroup of $G$ the measure preserving condition holds true because of the additional assumption in formula (\ref{hecke}).\end{proof}

We transform the formula in Proposition \ref{intersection} into a formula for the matrix coefficients for the action of $G$ on $L^{\infty}(\tilde{\Y}, \tilde{\nu})$.
We have, assuming the notations and definitions introduced above:

\begin{prop}
Let $\tilde{F}_0, \tilde{F}_1$ be measurable subsets of $\cZ_I$ (respectively $\cZ_{II}$).
Consider $F_0, F_1$ as in Proposition \ref{intersection}.
Then, using the action of $G$  introduced in the  above mentioned proposition, we have
$$
\tilde{\nu}(\tilde{F}_0 \cap g\tilde{F}_1) = \mathop{\lim}\limits_{n \to \infty} \mathop{\sum}\limits_{i,j} \nu_{\a}([(s_i^n)^{-1}[C_i^n \cap F_0]] \cap g[(s_j^n)^{-1}[C_j^n \cap F_1]]g^{-1}).
$$


\noindent On the right hand side we have  the action of $G$ on $\Gamma$. This is derived from the action of $G$ by conjugation on $\Gamma$ which extends to an action of $G$ on $^{\ast}\Gamma$ (partial action if $G \neq \Gamma$).

The formula remains valid for $F_0, F_1, \ldots, F_n$, subsets of $\tilde{\Y}$  and for group elements  $g_1, g_2, \ldots, g_n \in G$.
\end{prop}

\begin{proof}
This is very similar  the formula in Lemma \ref{intersection}. The only difference that needs a justification, is the fact that the action of $g \in G$ can be taken outside as it appears   in the right hand member of the equality.

In the case $G = \Gamma$, this is obvious, since we may substitute  $g\tilde{F}_0 $ by  $\tilde{gF_0}$, as $gF_0$ remains a $\Gamma$ - wandering domain (here the action of $g$ is derived from the action by conjugation on $\Gamma$ and $^{\ast}\Gamma$). Then the formula from Lemma \ref{intersection} gives  that the left hand side term is equal to
$$
\mathop{\lim}\limits_{n \to \infty}\mathop{\sum}\limits_{i,j}\nu_{\a}([(s_i^n)^{-1}[C_i^n \cap F_0]] \cap [(s_j^n)^{-1}[C_j^n \cap gF_0]]).
$$

But $g$ permutes the cosets $s_i^n \tilde{\Gamma}_n$, and hence  the action of $g$ may  be taken up front the parenthesis.

 In the case of type II, the argument works similarly, with the only difference that the  group   $\tilde{\Gamma}_n$ will be replaced by $\Gamma_n g\Gamma_0 g^{-1}$.
\end{proof}

\section{Model for the action of $G$ on $L^{\infty}(\tilde{\Y}, \tilde{\nu})$}

In this section we construct a model for the action of $G$ on $L^{\infty}(\tilde{\Y}, \tilde{\nu})$. We prove that the action may be realised in a space of the form $L^{\infty}(\Y_{\tilde{A}}, \nu_{\a})$, by changing the hyperfinite set $\cC_{\omega}(A_n)$ into a hyperfinite set $\cC_{\omega}(\tilde{A}_n)$ which sits in a suitable subset of the  fiber of $\pi_{A,K} : \Y_A \to K$, at $e$.

To do this note that the fiber of $\pi$ at $e$ is a reunion of subsets that  depend on the velocity which we impose on the ultrafilter convergence, for points in the fiber.

Choose a decreasing sequence of  subsets $U_n$  in $\omega$ with trivial intersection.
We introduce
$$
\pi_{\Gamma_n, U_n}^{-1}(e) = \{ (\gamma_n)_n \in \ ^{\ast}\Gamma \mid \gamma_n \in \Gamma_n \; {\rm if} \; n \in U_n \}.
$$

\noindent Here the subsets $(U_n)_{n\in \N}$ determine the "speed of convergence".

Consider the projection $\pi_{A,K} : \Y_A \to K$, and  consider the fiber of $\pi_{A,K}$ at a given $k\in K$.
Assume $k$ is represented as the intersection of the cosets $\overline{s_n\Gamma_n}, n\in \N$. The speed of convergence of some $(\gamma_n) \in \pi_{A,K}^{-1}(k)$ may be determined by the choice of the sets
$$
U_a = \{ n \mid \gamma_n \in s_a\Gamma_a \}, \; n \in \N.
$$

This gives a canonical definition for the "shape" of the convergence. We define ${\rm Aura}((\gamma_n)_{n\in \N})$ as the sequence
$$
\gamma_n^0 = s_a^{-1}\gamma_n \; {\rm if} \; n \in U_a \setminus U_{a + 1}.
$$
\noindent
Note that  $(\gamma_n^0)_{n\in \N}$ belongs to $\pi_{\Gamma_n, U_n}^{-1}(e)$.

We will do a similar construction for a hyperfinite set $\cC_{\omega}((A_n)_{n\in \N})$. Assume by eventually replacing with a smaller sequence of finite of subsets  $(A_n)_{n\in \N}$ that $\cC_{\omega}((A_n)_{n\in \N})$ is $\Gamma$ (respectively $\Gamma / \Gamma_0$ - wandering) for $\Gamma$ in $(\Y_A,\nu_\alpha)$.

 To do this one choses an enumeration $\gamma_1, \gamma_2,..., \gamma_a,\dots$ of  $\Gamma$ (respectively  of $\Gamma_0/\Gamma$, for $\Gamma_0\in \mathcal A_0$ in the type II case) and replaces for a sufficiently large sequence $(n_a)_{a \in \N}$ the set $A_n$
by 
$$A'_n= A_n\setminus \mathop\cup_{i=1,2, \dots a} \gamma_i A_n, n=n_a,\dots, n_{a+1}-1.$$
Then the measure space $(\Y_{A'},\nu_{\alpha})$ coincides with $(\Y_A,\nu_{\alpha})$ 

By the $\chi_1$ - saturation principle (\cite{Cut}), this amounts to an approximation, and hence it won't change the analysis of the states determining the continuity properties of  $\pi_{\rm Koop} \circ \Phi_{A, \omega}$.

\begin{defn}
Consider $\cC_{\omega}((A_n)_{n\in \N})$ as above.
Consider an enumeration $g_1, g_2, \ldots, g_a, \ldots$ of the group $G$. Fix $\varepsilon_a$ a sequence of positive numbers decreasing sufficiently fast to 0. Denote by $F$ the hyperfinite set  $\cC_{\omega}((A_n)_{n\in \N})$.

We let $U_a$ be the set of $n \in \N$ such that the sets $$({s_i^a})^{-1}(s_i^a\Gamma_a \cap A_n), \ i = 1, 2, \ldots, [\Gamma : \Gamma_a],$$ and their translates $$g_j\big[({s_i^a})^{-1}[s_i^a\Gamma_a \cap A_n]\big](g_j)^{-1}, j = 1, 2, \ldots, a$$ behave, together with their intersection, with respect to  the probability measure $$\nu_n = \frac{1}{\card A_n}(\mathop{\sum}\limits_{b \in A_n} \delta_b),$$ exactly, up to up to $\nu_n$ measure less than $\varepsilon_a$, as the sets $${(s_i^a)}^{-1}(C_i^a \cap F),$$
$$ g_j[{(s_i^a)}^{-1}(C_i^a \cap F)](g_j)^{-1}, j = 1, 2, \ldots, a,$$ and their intersections with respect to the measure $\nu_\alpha$. 

In particular up to $\nu_n$ measure less than  $\varepsilon_a$, the pieces ${(s_i^a)}^{-1}(s_i^a\Gamma_a \cap A_n)$ are disjoint, as they copy the behaviour of ${(s_i^a)}^{-1}(C_i^a \cap F)$.

Let  
\begin{equation}\label{aura}
{\rm Aura}(\cC_{\omega}((A_n)_{n\in \N}))=\cC_{\omega}((\tilde A_n)_{n\in \N}))
\end{equation}  be the hyperfinite set obtained from the following sequence of finite subsets of $\Gamma$:
$$
\tilde{A}_n = \mathop{\bigcup}\limits_{i = 1}^{[\Gamma : \Gamma_a]}(s_i^a)^{-1}[C_i^a \cap A_n] \subseteq \Gamma^a, \; n \in U_a \setminus U_{a + 1} .
$$
By the $\Gamma$-wandering property (respectively $\Gamma_0/\Gamma$- wandering property in the case of type II case)  we have 
$${\rm card} ((\tilde{A}_n)_{n \in \N})={\rm card\ } A\in ^\ast \N.$$

\end{defn}

\

\

Using the above definition, we obtain

\begin{prop} We use the above notations.
Let $F = \cC_{\omega}(A_n)$ and let  $\chi_{\tilde{F}}$ be the central support of $\chi_F$ in $\cZ = L^{\infty}(\tilde{\Y}, \tilde{\nu})$. Assume that $F$ is  a $\Gamma$ (respectively, for $\Gamma_0\in \mathcal A_0$, $\Gamma_0/\Gamma$) -wandering subset of $\Y_I$ (respectively $\Y_{\Gamma_0}$). In the second case we assume that  $F$ is  $\Gamma_0$-invariant. As noted above, by assuming this, we are not restricting the generality in the analysis of the Koopman representation of $C^{\ast}(G \rtimes (C^{\ast}(\Gamma \rtimes L^{\infty}(\Y_A, \nu_{\a}))))$.

Then the generalised moments
$$
\tilde{\nu}(\tilde{F} \cap g_1\tilde{F} \ldots g_n\tilde{F}), \; g_1, g_2, \ldots, g_n \in G, \; n \in \N
$$

\noindent are equal to the corresponding generalised moments associated to the  ${F}_0 = \cC_{\omega}(\tilde A_n)$, associated to  ${\rm Aura}(\cC_{\omega}((A_n)_{n\in \N}))$ with respect to the Loeb measure $\nu_{{\rm card\ }\tilde{A}}$.

In particular the action of $G$ on $\cZ = L^{\infty}(\tilde{\Y}, \tilde{\nu})$ is amenable, if and only if the action of $G$ on $L^{\infty}(\Y_{\tilde{A}}, \nu_{{\rm card\ }\tilde{A}})$ is amenable. In this statement by amenability we understand nuclearity of the corresponding crossed product algebras.
\end{prop}

\begin{proof}
Indeed choosing $(\varepsilon_a)_{a \in \N}$ a  sufficiently fast decreasing sequence, this is a consequence of the relations imposed in the choice of the hyperfinite set introduced in formula (\ref{aura}).
\end{proof}

We assume in the rest of this section  that  $G = \Gamma$, that $\Gamma$ is exact and it is an  i.c.c. group.

\begin{proof}[Proof of Theorem \ref {AOconditions}]
Indeed, since $G$ acts by conjugation on $\pi^{-1}(e)$, if we choose $\Gamma_n$, so that the points in $\Gamma_n$ have amenable stabilizers then the action of $\Gamma$ on $(\Y_{\tilde{A}}, \nu_{\a})$ is amenable, since $\Y_{\tilde{A}} \subseteq \pi_{\Gamma_n, U_n}^{-1}(e)$.

This implies that $\Gamma$ acts amenable on the center of the nuclear algebra\break $\{ \pi_{\rm Koop}(C^{\ast}(\Gamma \rtimes C^{\ast}(\Gamma))) \}''$ and hence, by (\cite{CAD}, Proposition 8.2), the $C^{\ast}$ - algebra $$\{ \pi_{\rm Koop}(C^{\ast}(\Gamma \rtimes C^{\ast}(\Gamma \rtimes L^{\infty}(y_A, \nu_A)))) \}''$$ is nuclear. Inside we have the image of the crossed product algebra
$$
(\pi_{\rm Koop} \circ \Phi_{A, \omega})(C^{\ast}(\Gamma \rtimes C^{\ast}(\Gamma \rtimes l^{\infty}(\Gamma)))).
$$

By the nuclearity of the larger crossed product, it follows that the representation of the smaller crossed product algebra factorizes to the reduced $C^{\ast}$ - algebra
$$
C^{\ast}_{\rm red}(\Gamma \rtimes C^{\ast}(\Gamma \rtimes l^{\infty}(\Gamma))) = C^{\ast}_{\rm red}(\Gamma \rtimes C^{\ast}_{\rm red}(\Gamma \rtimes l^{\infty}(\Gamma))) = 
$$
$$
= C^{\ast}_{\rm red}((\Gamma \rtimes \Gamma) \rtimes l^{\infty}(\Gamma)).
$$

\end{proof}



\section{The non free part of the  action of $(G \rtimes \Gamma)$ on 
$^\ast \Gamma$}\label{nonfree}

In this section we treat the possible cases when the action of $G \rtimes \Gamma$ is non free. We will consider only the cases $G \rtimes \Gamma = \Gamma \times \Gamma^{\rm op}$ for $\Gamma = \PGL_2(\Z[\frac{1}{p}])$ or $\Gamma = \SL_3(\Z)$.

Consider the subset $(\Gamma^{\ast})^{\rm fixed}$ of $^{\ast}\Gamma$ consisting of all points in $^{\ast}\Gamma$ having a nontrivial stabilizer in $\Gamma\times \Gamma^{\rm op}$. For $(\Y_A, \nu_{\a})$  as in Definition  \ref {$y_a$} , let $$\Y_A^{\rm fixed} = \Y_A \cap (^{\ast}\Gamma)^{\rm fixed}.$$
\noindent
Let $\Phi_{A} : C(\partial(\beta\Gamma))\to \B(L^2((\Y_A, \nu_{\a})$ be the map constructed in formula (\ref{phir}) in Section \ref{intro}.
The main result of this section is the fact that measure spaces of the form $\Y_A^{\rm fixed}$ and hence the set $(^{\ast}\Gamma)^{\rm fixed}$   are contained in the image, through the morphism $\Phi_{A}$, of a canonical subset of $\ell^\infty(\Gamma)$.
The main statement of this section is:  

\begin{lemma}\label{y1}
Let $\partial(\beta\Gamma)^{\rm fixed}$ be the subset of $\partial(\beta\Gamma)$ obtained as a reunion of   cosets of the form $[\Gamma_1 x]$, where $\Gamma_1$ is a centralizer group as in Lemma \ref{fixed} below.

Then $\Y_A^{\rm fixed}$ is equal to the image of $\Phi_A(\chi_{\partial(\beta\Gamma)^{\rm fixed}})$.
\end{lemma}


\noindent Before proving the lemma we recall a few  standard facts.
%
%
For $C \subseteq \Gamma$ we denote by $C'$ the subgroup of all $\gamma$ in $\Gamma$ such that $\gamma c \gamma^{-1} = c$, for all $ c \in C$.
If $C$ is a subgroup, then $C'$ is the centralizer subgroup.

\begin{lemma}\label{fixed}
Let $\Gamma \times \Gamma^{\rm op}$ act on $\Gamma$ by left and right multiplication.
Fix $(\gamma_1, \gamma_2) \in \Gamma \times \Gamma^{\rm op}$.
Assume that $(\gamma_1, \gamma_2)$ keeps  $x$ fixed, that is $\gamma_1 x \gamma_2^{-1} = x$ (equivalently $\gamma_2 = x^{-1}\gamma_1 x$).

Then the set of points fixed by $\{ \gamma_1, \gamma_2 \}$ is the coset $\{ \gamma_1 \}'x = x \{ \gamma_2 \}'$.
\end{lemma}

\begin{proof}
Assume $y$ is another point fixed by $(\gamma_1, \gamma_2)$. Then $\gamma_2 = x^{-1}\gamma_1 x = y^{-1}\gamma_1 y$ and hence $yx^{-1}$ commutes to $\Gamma_1=\{ \gamma_1 \}'$ and hence $y$ belongs to $\{ \gamma_1 \}'x = x\{ \gamma_2 \}'$.

The reciprocal is clear.
\end{proof}

The following lemma is known statement about intersections of cosets.

\begin{lemma}
Let $H$ be a discrete group. Let $H_0, H_1$ be two subgroups. We consider two cosets $H_0x_0, H_1x_1$, $x_0, x_1 \in H$.

Assume the intersection $H_0x_0 \cap H_1x_1$ is non-void, and let $x$ be a point in the intersection.

Then $H_0x_1 \cap H_1x_1 = (H_0 \cap H_1)x$.
\end{lemma}

\begin{proof}
The term on the right hand side is obviously contained in the term of the left hand side of the equality.
We prove that the intersection of the right hand side is contained in the term $(H_0 \cap H_1)x$.

Let $x' \neq x$ be any other element of the intersection $H_0x_0 \cap H_1x_1$ which by hypothesis contains $x$.

Then there exist $h_0, h_0' \in H_0$ (respectively $h_1, h_1' \in H_1$) such that
$$
\begin{array}{cc}
x = h_0x_0 = h_1x_1\\[3mm]
x' = h_0'x_0 = h_1'x_1
\end{array}
$$

Then
$$
\begin{array}{cc}
x(x')^{-1} = (h_0x_0)(h_0'x_0') = h_0(h_0')^{-1} \in H_0\\[3mm]
x(x')^{-1} = (h_1x_1)(h_1'x_1)^{-1} = h_1(h_1')^{-1} \in H_1
\end{array}
$$

Hence $\T = h_0(h_0')^{-1} = h_1(h_1')^{-1}$ to $H_0 \cap H_1$.

Then
$$
x(x')^{-1} = \T \in H_0 \cap H_1
$$

\noindent and hence $x' = \T^{-1}x$ which consequently belongs to $(H_0 \cap H_1)x$. 
\end{proof}

We analyse the action of the group $\Gamma \times \Gamma^{\rm op}$ on the characteristic functions of cosets as above, in the  algebra  $C(\partial(\beta\Gamma))$.

Because of the previous lemma, if $\Gamma_0, \Gamma_1$ are subgroups of $\Gamma$ with trivial intersection, then $\Gamma_0x \cap \Gamma_1y$ consists of at most one point for all $x, y \in \Gamma$.

\begin{proof}[Proof of Lemma \ref {y1}]
This is a consequence of Lemma \ref{fixed}. Indeed if $\gamma= (\gamma_n) \in \ ^{\ast}\Gamma$ is left invariant by $(\gamma_1, \gamma_2)$, then $\gamma_1 a_n \gamma_2^{-1} = a_n$ for infinitely many $n$, and hence $(a_n)$ belongs eventually to the image through $\Phi_A$ of one of the coset $[\Gamma_1 x]$ (in fact in $\pi_{\Y_A}^{-1}([\Gamma_1 x])$) where $\pi_{\Y_A} : \Y_A^{\rm fixed} \to \partial(\beta\Gamma)^{\rm fixed}$ is the canonical projection.
\end{proof}

\section{the group $\PGL_2(\Z[\frac{1}{p}])$}

In the case $\Gamma = \PGL_2(\Z[\frac{1}{p}])$  we prove that the conditions of Theorem \ref{AOconditions} hold true.  Hence the Koopmann representation factorizes to $C^{\ast}_{\rm red}(\Gamma \times \Gamma^{\rm op})$.

\begin{thm}
$\Gamma = \PGL_2(\Z[\frac{1}{p}])$ has the ${\rm AO}$ property.
\end{thm}

\begin{proof}
In this case there are no points with non-amenable stabilizers.  Using  Theorem \ref{AOconditions} in previous section, to it follows that is sufficient  analyze the action of $\Gamma \times \Gamma^{\rm op}$ on $\Y_A^{\rm fixed}$.

The possible groups $\Gamma_1' = \{ g \}'$, $g \in \Gamma$, as considered  in Lemma \ref{fixed} are determined as follows: 



Fix $g$ an element of $\PGL_2(\Z[\frac1{p}])$. 
There are two cases: either $g$ viewed
as a matrix with real entries has two distinct eigenvalues, or either $g$ is conjugated 
to an element in the triangular group
$$
T_p = \left\{ \left( \begin{array}{cc} a & b \\ 0 & a \end{array}
\right) \mid a,b \in \Z[\tfrac1{p}]\right\},
$$
considered as subgroup of $ \PGL(2, \Z[\tfrac1{p}])$.
Note that here, since we are using the projective groups, 
all matrices are  considered, as classes modulo the scalar matrices.

In the first case, the commutant of $g$ will be either finite (e.g., if $g$ is conjugate
to $\left(\begin{array}{cc} 0 & -1 \\ 1 & 0 \end{array}
\right)$) or a maximal abelian subgroup of $\PGL_2(\Z[\frac1{p}])$ with 
trivial normaliser (and hence isomorphic to $\Z$).

In the second case the commutant will be the group $T_p$ itself. It is obvious to see that
$T_p$ is a maximal abelian group with trivial normaliser.

The group 
 ${\Gamma_1}$ defined in Lemma \ref{fixed} is either of the form 
$$
\begin{array}{ll}
\Gamma'_1 = \{g^n\},\ \hbox{\ if \ }  g\in \Gamma \hbox{ has distinct eigenvalues and}\\
\Gamma'_1 \cong \Z, \ \Gamma_1 \hbox{ maximal abelian} 
\end{array}\leqno(\alpha)
$$
or either
$$
\hbox{ the group $\Gamma'_1$ is a conjugate of $T_p$.} \leqno(\beta)
$$

Clearly, for two subgroups as in property $(\alpha)$,  since they are infinite maximal abelian, if they have infinite intersection, than they coincide. No group of the type in $(\alpha)$
can intersect (except in the trivial element) a group in $(\beta)$.

A simple computation shows that if $g$ belongs to $\Gamma = \PGL_2(\Z[\frac1{p}])$
and $g T_p g^{-1} \cap T_p$ is non-trivial, then $g$ must belong to $T_p$ (this is
a stronger property than having trivial normaliser). 


Thus $\Gamma \times \Gamma^{\rm op}$ on $\partial(\beta\Gamma)^{\rm fixed}$ acts as follows. Fix $\Gamma_1'$ as above.

Then the set $$\gamma_1[\Gamma_1' x]\gamma_2^{-1}\subseteq \partial(\beta\Gamma)^{\rm fixed}, $$ will be $\Gamma \times \Gamma / \Gamma_1' \times x^{-1}\Gamma_1' x$ wandering. Moreover $\Gamma_1' \times x^{-1}\Gamma_1' x$ acts amenably on the set $[\Gamma_1' x]$ as $\Gamma_1'$ is amenable. Thus the Koopmann representation is weakly contained in the quasiregular representation of $\Gamma \times \Gamma^{\rm op}$ on $\Gamma \times \Gamma^{\rm op}$ on $l^2(\Gamma \times \Gamma^{\rm op}/ \Gamma_1' \times x\Gamma_1' x^{-1})$.

The latest, since $\Gamma_1'$ is amenable is contained in the left regular representation of $\Gamma \times \Gamma^{\rm op}$.

\end{proof}

\begin{cor}\label{nonao} 
The group $\Gamma = \PGL_2(\Z[\frac1{p}])$ has  the property AO but does not have the property
$\mathcal S$ of Ozawa.
\end{cor}

\begin{proof}
As Sergey Neshveyev and Makoto Yamashita kindly pointed out to us, the group $\Gamma$ does not have
the stronger related property $\mathcal S$ of Ozawa. Indeed, being a lattice 
in $\PSL_2(\R)\times \PGL_2(\mathbb Q_p)$ (because of  [Ih]), it is stably measurably equivalent
to $F_2 \times F_2$. But as proven by Sako [Sa], the property $\mathcal S$
is preserved by stably measurable equivalence, and since $F_2 \times F_2$ does not
have this property, it follows that $\PGL_2(\Z[\frac1{p}])$ does not have property $\mathcal S$ of Ozawa ([Oz]),
but does have the property AO.
\end{proof}

\vspace{0.2cm}



\section{Examples: the case of $\SL_3(\Z)$}

We will adapt the conditions of Lemma \ref{fixed} for the group $\SL_3(\Z)$. 
For this purpose we introduce the following subgroups of $\SL_3(\Z)$.

Let $H$ be the Heisenberg subgroup consisting of all matrices of the form
$$
\left( \begin{array}{ccc} 1 & * & * \\ 0 & 1 & * \\ 0 & 0 & 1 \end{array} 
\right),
$$
with integer entries. Let $\SL_2(\Z) \subseteq \SL_3(\Z)$ be the 
canonical representation of $\SL_2(\Z)$ as a subgroup of $\SL_3(\Z)$. Thus $\SL_2(\Z)$  is
the set of all matrices in $\SL_3(\Z)$ of the form
$$
\left( \begin{array}{ccc} * & * & 0 \\ * & * & 0 \\ 0 & 0 & 1 \end{array} 
\right).
$$

Let $E$ be the matrix $$
\left( \begin{array}{ccc} -1 & 0 & 0 \\ 0 & -1 & 0 \\ 0 & 0 & 1 \end{array} 
\right)
.$$ Let $H_2 = H\cap \SL_2(\Z)$. This is the abelian subgroup of triangular matrices.

As for $H_2$, the subgroup $H$ has the property
that for $\gamma$ in $\SL_3(\Z)\setminus H$ the intersection
$\gamma H \gamma^{-1} \cap H $ is the trivial subgroup.


We  analyze the action of $\Gamma \times \Gamma^{\rm op}$ on $(\Y_A^{\rm fixed}, \nu_{\a})$.
The only element having non-amenable stabilizers (along with all elements in its conjugacy orbit) is $E$. But the modular subgroups separate the orbit of $E$ from the identity.

To prove that  $\SL_3(\Z)$ has the AO property it remains to  analyze the action of $\Gamma$ on $\Y_A^{\rm fixed}$.
%
In the case of $\SL_3(\Z)$, differently from the case of $\PGL_2(\Z[\frac1{p}])$,
the commutant of $E$ is equal to $\SL_2(\Z)$, a non-amenable group, and moreover the
intersections $g^{-1}\SL_2(\Z)g \cap \SL_2(\Z)$ might be non-trivial,
and  infinite, for $g$ not belonging to $\SL_2(\Z)$. However, as we prove in the next statement,
each of the above intersections  will be a subgroup   of a conjugate of the group $H$.
More precisely, we have:

\begin{lemma}\label{Tausky}
Assume that $x,y$ are non-trivial elements of $\SL_2(\Z) \subseteq \SL_3(\Z)$
and $g$ belongs to $\SL_3(\Z) \setminus \SL_2(\Z)$ such that
$g x g^{-1} = y$.

Then  there exists $\gamma$ in $\SL_2(\Z)$ such that $\gamma^{-1}x\gamma$
belongs to $H_2$ and there exists $\gamma_0$ in $\SL_2(\Z)$ and $h$ in $H$
such that $g = \gamma_0 (\gamma h \gamma^{-1})$.

In particular any subgroup obtained as  non trivial intersection of 
$$ \SL_2(\Z)\cap g(\SL_2(\Z))g^{-1}, g \in \SL_3(\Z) \setminus \SL_2(\Z),$$
is contained in a subgroup 
$\tilde {H}_2$,  which is a conjugate in $\SL_2(\Z)$ to the group $H_2$.
\end{lemma}

\begin{proof}

By the results of Olga Tausky (\cite{OT},  see also \cite{LMD} and  the references in there), the
conjugacy classes for elements in $\SL_3(\Z)$ are determined by ideal classes
in the ring obtained by adjoining to $\Z$ the roots of the characteristic polynomial.

Hence if $x,y$ belong to $\SL_2(\Z)$ and are conjugated in $\SL_3(\Z)$, they are
also conjugate in $\SL_2(\Z)$ and hence there exists $\gamma_0$ in $\SL_2(\Z)$ 
such that 
$$
g x g^{-1} = \gamma_0 x \gamma_0^{-1} = y.
$$
But then $(\gamma_0^{-1} g)x (\gamma_0^{-1} g)^{-1} = x$ and hence
$\gamma_0^{-1}g\in \SL_3(\Z) \setminus \SL_2(\Z)$ commutes with $x$.

The only cases, when an element in $\SL_2(\Z)$ has something in the commutant, which belongs to $\SL_3(\Z) \setminus \SL_2(\Z)$, are the conjugates by elements in $\SL_2(\Z)$ of the following example:

\vspace{0.2cm}

The group $H_2$ commutes with the group
$$
H_2^0 =
\left\{ 
\left(
\begin{array}{ccc} 
1 & 0 & 0\\[3mm]
0 & 1 & n\\[3mm]
0 & 0 & 1
\end{array}
\right) \mid n \in \Z 
\right\}.
$$

Hence we may assume $x = \a h_2 \a^{-1}$, $h_2 \in H_2$, $\a \in \SL_2(\Z)$
$$
g^{-1}\gamma_0 = \a()h_2^0)^{-1}\a^{-1}, \; h_2^0 \in H_2^0.
$$

Thus $g = \gamma_0\a(h_2^0)\a^{-1}$.

But in this case
$$
\PSL_2(\Z) \cap g\PSL_2(\Z)g^{-1} = 
$$
$$
= \gamma_0 \a [\SL_2(\Z) \cap (h_2^0)\a^{-1} \PSL_2(\Z) \a(h_2^0)^{-1}]\a^{-1}\gamma_0^{-1} =
$$
$$
= \gamma_0 \a [\SL_2(\Z) \cap (h_2^0)\PSL_2(\Z)h_2^0]\a^{-1}\gamma_0^{-1} =
$$
$$
= \gamma_0\a [H_2]\a^{-1}\gamma_0^{-1}.
$$

Since $\gamma_0 \a \in \SL_2(\Z)$ it follows that $\PSL_2(\Z) \cap g\PSL_2(\R)g^{-1}$ is equal to a conjugate by an element in $\SL_2(\Z)$ of $H_2$

\end{proof}

\begin{thm} The group 
$\SL_3(\Z)$ has the ${\rm AO}$ property.
\end{thm}

\begin{proof}
We analyze $\{g\}'$ if $g \in \SL_3(\Z)$.
If $g$ has three distinct eigenvalues the situation is exactly as  in the case of $\PSL_2(\Z[\frac{1}{p}])$.

If $g$ is conjugated to an element in $H_3$ the commutant is a possible larger subgroup of $H_3$ (as the commutant of $H_2$ contains $H_2$ and $H_2^0$).
Since $H_3$ is amenable, we reapply the argument from $\PSL_2(\Z[\frac{1}{p}])$ for the subset of $\partial(\beta\Gamma)$ generated by subsets of the form $[H_3 x]$.
We use the fact that  $H_3$ is amenable, and $H_3 \cap g H_3 g^{-1}$ is trivial unless $g \in H_3$.  The argument from the case of $\PSL_2(\Z[\frac{1}{p}])$ works again in this case.

The remaining situation is when $(\gamma_1, \gamma_2)$ is of the form $(E, x^{-1} E x)$ fixing $x$.
In this case the points fixed by such an element are $$[\SL_2(\Z)x] =[x(x^{-1}\SL_2(\Z  )x)], x \in \SL_3(\Z).$$

We intersect the subset of $\Y_A^{\rm fixed}$ with the complement of the subset generated by the cosets $[gH_3g^{-1} x], x, g\in \SL_3(\Z) $.  Because of Lemma \ref{Tausky} it follows that in the above set difference, the set
$[\PSL_2(\Z)x]$ is $$\Gamma \times \Gamma^{\rm op} / \SL_2(\Z) \times x^{-1}\SL_2(\Z)^{\rm op}x$$ wandering. Moreover the action of $\SL_2(\Z) \times x^{-1}\SL_2(\Z)^{\rm op}x$ on $\PSL_2(\Z)x$ is equivalent to  action of $\SL_2(\Z) \times \SL_2(\Z)^{\rm op}$ on $\SL_2(\Z)$.

Because AO for $\SL_2(\Z)$ it follows that $\SL_3(\Z) \times \SL_3(\Z)^{\rm op}$ acts amenably  on  $(\Y_A)^{\rm fixed}$.
\end{proof}

We are indebted to Kang Li for pointing us out the following corollary  of the  property in the previous theorem. To obtain the result one uses
  the results in  Skandalis's   paper on the AO property (\cite {Sk}, see also \cite {OU}):
\begin{cor*}The  full group  C$^\ast$-algebra  C$^\ast(\SL_3(\Z))$ is not K-exact.
\end{cor*}

{\bf Acknowledgement.} This article was written during the stay of the author in April-May, 2011 at the IHES Institute, Bures sur Yvette, and at  the IHP Institute, Paris, May-June 2011, during the Workshop ``Von Neumann algebras and ergodic theory of group actions". The author is grateful for the warm welcome at both institutes. The author is indebted to Maxim Kontsevich,  Thibault Damour and
Laurent Lafforgue for discussions, around the subject of this paper, during his stay at IHES. The author is indebted to Ionu\c t Chifan,
Adrian Ioana, Narutaka Ozawa, Sergey Neshveyev, Stefaan Vaes Marius Junge, Nicolas Monod, Peter Loeb, Uffe Haagerup, Ryszard Nest, Alexander Gorodnik, Peter Loeb, Grigore Ciurea, Ovidiu Pasarescu, Kang Li for comments on  the present  article.

\end{document}